\documentclass[11pt,a4paper,leqno,twoside]{amsart}
\usepackage{latexsym,amssymb,amsmath}
\input xy
\xyoption{all}

\def\CC{\mathbb C}

\def\RR{\mathbb R}
\def\HH{\mathbb H}
\def\AA{{\mathbb A}}

\def\OO{\mathbb O}

\def\11{\mathbf 1}
\def\PP{\mathbb P}

\def\e1{\varepsilon_1}
\def\e2{\varepsilon_2}
\def\e3{\varepsilon_3}

\def\P2{{\PP}^2}

\def\00{\underline{0}}
\def\J0{{\cal J}_3(\underline{0})}

\def\PJ0{\PP({\cal J}_3(\underline{0}))}

\def\e{\varepsilon}

\def\AP2{{\AA\PP}^2}
\def\RP2{{\RR\PP}^2}
\def\CP2{{\CC\PP}^2}
\def\HP2{{\HH\PP}^2}
\def\OP2{{\OO\PP}^2}

\newtheorem{theo}{Theorem}[section]

\newtheorem{lemm}[theo]{Lemma}
\newtheorem{prop}[theo]{Proposition}
\newtheorem{conj}[theo]{Conjecture}

\newtheorem{fact}[theo]{Fact}
\newtheorem{prob}[theo]{Problem}

\theoremstyle{remark}
\newtheorem{rema}[theo]{Remark}

\begin{document}
\title[Brauer $p$-dimension of HDV-fields of residual chacteristic $p$]{On the Brauer $p$-dimension of Henselian discrete valued
fields of residual characteristic $p > 0$}

\author{Ivan D. Chipchakov}
\address{Institute of Mathematics and Informatics\\Bulgarian Academy
of Sciences\\1113 Sofia, Bulgaria: E-mail address:
chipchak@math.bas.bg}

\begin{abstract}
Let $(K, v)$ be a Henselian discrete valued field with residue field
$\widehat K$ of characteristic $p > 0$, and Brd$_{p}(K)$ be the Brauer
$p$-dimension of $K$. This paper shows that Brd$_{p}(K) \ge n$ if
$[\widehat K\colon \widehat K ^{p}] = p ^{n}$, for some $n \in
\mathbb{N}$. It proves that Brd$_{p}(K) = \infty $ if and only if
$[\widehat K\colon \widehat K ^{p}] = \infty $.
\end{abstract}

\maketitle

\par
\medskip
\section{\bf Introduction}
\par
\medskip
Let $E$ be a field, Br$(E)$ its Brauer group, $s(E)$ the class of
associative finite-dimensional central simple algebras over $E$, and
$d(E)$ the subclass of division algebras $D \in s(E)$. For each
$A \in s(E)$, let $[A]$ be the equivalence class of $A$ in Br$(E)$,
and let deg$(A)$, ind$(A)$, exp$(A)$ be the degree, the Schur
index and the exponent of $A$, respectively. It is well-known (cf.
\cite{P}, Sect. 14.4) that exp$(A)$ divides ind$(A)$ and shares with
it the same set of prime divisors; also, ind$(A) \mid {\rm deg}(A)$,
and deg$(A) = {\rm ind}(A)$ if and only if $A \in d(E)$. Note that if
$B _{1}, B _{2} \in s(E)$ and g.c.d.$\{{\rm ind}(B _{1}), {\rm ind}(B
_{2})\} = 1$, then ind$(B _{1} \otimes _{E} B _{2}) = {\rm ind}(B
_{1}){\rm ind}(B _{2})$; equivalently, if $B _{j} ^{\prime } \in
d(E)$, $j = 1, 2$, and g.c.d.$\{{\rm deg}(B _{1} ^{\prime }), {\rm
deg}(B _{2} ^{\prime })\}$ $= 1$, then $B _{1} ^{\prime } \otimes
_{E} B _{2} ^{\prime } \in d(E)$ (see \cite{P}, Sect. 13.4). Since
Br$(E)$ is an abelian torsion group and ind$(A)$, exp$(A)$ are
invariants both of $A$ and $[A]$, these results show that the study
of the restrictions on the pairs ind$(A)$, exp$(A)$, $A \in s(E)$,
reduces to the special case of $p$-primary pairs, for an arbitrary
prime $p$. The Brauer $p$-dimensions Brd$_{p}(E)$, $p \in
\mathbb P$, where $\mathbb P$ is the set of prime numbers, are
defined as in \cite{ABGV}, and contain essential information on these
restrictions. We say that Brd$_{p}(E) = n < \infty $, for a given $p
\in \mathbb P$, if $n$ is the least integer $\ge 0$, for which ind$(P)
\mid {\rm exp}(P) ^{n}$ whenever $P \in s(E)$ and $[P]$ lies in the
$p$-component Br$(E) _{p}$ of Br$(E)$; if no such $n$ exists, we put
Brd$_{p}(E) = \infty $. For instance, Brd$_{p}(E) \le 1$, for all $p
\in \mathbb P$, if and only if $E$ is a stable field, i.e. deg$(D) =
{\rm exp}(D)$, for each $D \in d(E)$; Brd$_{p'}(E) = 0$, for some $p
^{\prime } \in \mathbb P$, if and only if the $p'$-component Br$(E)
_{p'}$ of Br$(E)$ is trivial.
\par
The absolute Brauer $p$-dimension abrd$_{p}(E)$ of $E$ is defined
to be the supremum of Brd$_{p}(R)\colon R \in {\rm Fe}(E)$, where
Fe$(E)$ is the set of finite extensions of $E$ in a separable
closure $E _{\rm sep}$. This trivially implies abrd$_{p}(E) \ge
{\rm Brd}_{p}(E)$, for each $p$. We have abrd$_{p}(E) \le 1$, $p
\in \mathbb P$, if $E$ is an absolutely stable field, i.e. its
finite extensions are stable fields. Class field theory gives
examples of such fields: it shows that Brd$_{p}(\Phi ) = {\rm
abrd}_{p}(\Phi ) = 1$, $p \in \mathbb P$, if $\Phi $ is a global
or local field (see, e.g., \cite{Re}, (31.4) and (32.19)). The
same equalities hold, if $\Phi = \Phi _{0}((X))((Y))$ is an
iterated formal Laurent power series field in $2$ variables over a
quasifinite field $\Phi _{0}$ (see \cite{Ch1}, Corollary~4.5
(ii)).
\par
The knowledge of the sequence Brd$_{p}(E), {\rm abrd}_{p}(E)\colon p
\in \mathbb P$, is helpful for better understanding the behaviour of
index-exponent relations over finitely-generated transcendental
extensions of $E$ \cite{Ch4}. This is demonstrated by the description
in \cite{Ch5} of the set of sequences Brd$_{p}(K _{q})$, abrd$_{p}(K
_{q})$, $p \in \mathbb P$, $p \neq q$, where $K _{q}$ runs across the
class of fields with Henselian valuations $v _{q}$ whose residue
fields $\widehat K _{q}$ are perfect of characteristic $q \ge 0$,
such that their absolute Galois groups $\mathcal{G}_{\widehat K _{q}}
= \mathcal{G}(\widehat K _{q,{\rm sep}}/\widehat K _{q})$ are
projective profinite groups, in the sense of \cite{S1}. The
description relies on formulae for Brd$_{p}(K _{q})$, $p \neq q$,
which depend only on whether $\widehat K _{q}$ contains a primitive
$p$-th root of unity. Thus Brd$_{p}(K _{q})$ is determined, for each
$p \neq q$, by two invariants: one of the value group $v _{q}(K
_{q})$, and one of the Galois group $\mathcal{G}(\widehat
K _{q}(p)/\widehat K _{q})$ of the maximal $p$-extension $\widehat
K _{q}(p)$ of $\widehat K _{q}$ in $\widehat K _{q,{\rm sep}}$.
\par
A formula for Brd$_{q}(K _{q})$ in terms of invariants of $\widehat
K _{q}$ and $v(K _{q})$ has also been found when char$(K _{q}) = q >
0$, $\widehat K _{q}$ is perfect and $(K _{q}, v _{q})$ is a
maximally complete field (see \cite{Ch6}, Proposition~3.5). By
definition, the imposed restriction on $(K _{q}, v _{q})$ means that
it does not admit immediate proper extensions, i.e. valued extensions
$(K _{q} ^{\prime }, v'_{q}) \neq (K _{q}, v _{q})$ with $\widehat K
_{q} ^{\prime } = \widehat K _{q}$ and $v'_{q}(K _{q} ^{\prime }) = v
_{q}(K _{q})$. The considered fields are singled out by the fact
(established by Krull, see \cite{Wa}, Theorem~31.24 and page 483)
that every valued field $(L _{0}, \lambda _{0})$ has an immediate
extension $(L _{1}, \lambda _{1})$ that is a maximally complete
field. Note here that no formula for Brd$_{q}(K _{q})$ as above
exists if $(K _{q}, v _{q})$ is only Henselian. More precisely, one
can show using suitably chosen valued subfields of maximally complete
fields that if $(K, v)$ runs across the class of Henselian fields of
characteristic $q$, then Brd$_{q}(K)$ does not depend only on
$\widehat K$ and $v(K)$. Specifically, it has been proved (see
\cite{Ch6}, Example~3.7) that for any integer $t \ge 2$, the iterated
formal Laurent power series field $Y _{t} = \mathbb{F} _{q}((T _{1}))
\dots ((T _{t}))$ in $t$ variables over the field $\mathbb{F} _{q}$
with $q$ elements possesses subfields $K _{\infty }$ and $K _{n}$, $n
\in \mathbb{N}$, such that:
\par
\medskip\noindent
(1.1) (a) Brd$_{q}(K _{\infty }) = \infty $; $n + t - 1 \le {\rm
Brd}_{q}(K _{n}) \le n + t$, for each $n \in \mathbb{N}$;
\par
(b) The valuations $v _{m}$ of $K _{m}$, $m \le \infty $, induced
by the standard $\mathbb{Z} ^{t}$-valued valuation of $Y _{t}$ are
Henselian with $\widehat K _{m} = \mathbb{F} _{q}$ and $v _{m}(K
_{m}) = \mathbb{Z} ^{t}$; here $\mathbb{Z} ^{t}$ is viewed as an
abelian group endowed with the inverse-lexicographic ordering.
\par
\medskip
Statement (1.1) attracts interest in the study of Brauer
$p$-dimensions of Henselian fields of residual characteristic $p > 0$
from suitably chosen special classes. This paper considers
Brd$_{p}(K)$, for a Henselian discrete valued field (abbr., an
HDV-field) $(K, v)$ with char$(\widehat K) = p$. Our research is
related to the problem of describing index-exponent relations over
finitely-generated field extensions. It proves the right-to-left
implication in the equivalence \par\noindent Brd$_{p}(K) = \infty
$ $\Leftrightarrow $ the degree $[\widehat K\colon \widehat K
^{p}]$ is infinite (in case char$(K) = 0$, the inverse implication is
a consequence of \cite{PS}, Corollary~2.5, see also Fact \ref{fact3.5}),
$\widehat K ^{p}$ being the subfield $\{u ^{p}\colon u \in \widehat
K\}$. When $[\widehat K\colon \widehat K ^{p}] < \infty $, we prove
the lower bound in the following conjecture (stated by Bhaskhar and
Haase \cite{BH}\footnote{The Brauer $p$-dimension, in the sense of
\cite{PS} and \cite{BH}, means the same as the absolute Brauer
$p$-dimension in the present paper.} for complete discrete valued
fields):
\par
\smallskip
\begin{conj}
\label{conj1.1} If $(K, v)$ is an {\rm HDV}-field with {\rm
char}$(\widehat K) = p > 0$ and $[\widehat K\colon \widehat K
^{p}] = p ^{n}$, for some $n \in \mathbb{N}$, then $n \le {\rm
abrd}_{p}(K) \le n + 1$.
\end{conj}
\par
\medskip\noindent
Conjecture \ref{conj1.1} has been stated at the end of \cite{BH},
under the extra hypothesis that char$(K) = 0$ and char$(\widehat
K) = p$. This restriction is not emphasized in the present paper,
as we prove Conjecture \ref{conj1.1} in case char$(K) = p$ (see
Proposition \ref{prop7.1}). Note also that the class of HDV-fields
$(K, v)$ with char$(\widehat K) = p > 0$ and $[\widehat K\colon
\widehat K ^{p}] = p ^{n}$ is closed under taking finite
extensions (cf. \cite{E3}, Corollary~14.2.2, and \cite{BH},
Lemma~2.12). It is therefore clear that the upper bound in
Conjecture \ref{conj1.1} will follow, if the inequality
Brd$_{p}(K) \le n + 1$ holds, for an arbitrary HDV-field $(K, v)$
with $\widehat K$ as above. The inequality Brd$_{p}(K) \ge n$
implies trivially the lower bound $n \le {\rm abrd}_{p}(K)$ in
Conjecture \ref{conj1.1}. It attracts interest in finding formulae
for Brd$_{p}(K)$, for example, when $(K, v)$ belongs to some basic
classes of HDV-fields of residual characteristic $p$ (see
Conjecture \ref{conj7.3} and Problem \ref{prob7.4}).

\par\vskip0.8truecm\noindent

{\bf Basic notation and abbreviations used in the paper}

\begin{itemize}

\item $\mathbb{P}$ - the set of prime numbers; $\mathbb{N}$ - the
set of positive integers; $\mathbb{Z}$ - the set (additive group,
ring) of integers;

\item $\mathbb{Q}$ and $\mathbb{R}$ the additive groups (the fields)
of rational numbers and of real numbers, respectively;

\item Abbreviations: HDV - Henselian discrete valued; TR - totally
ramified;

\item For any field $E$, we use the following notation:

\item $E ^{\ast }$ is the multiplicative group of $E$; $E ^{\ast n}$
is the subgroup of $n$-th powers $E ^{\ast n} = \{\alpha
^{n}\colon \alpha \in E ^{\ast }\}$, for each $n \in \mathbb{N}$;

\item $s(E)$ - the class of associative finite-dimensional central
simple $E$-algebras, $d(E)$ - the subclass of division algebras $D
\in s(E)$, Br$(E)$ - the Brauer group of $E$;

\item $E _{\rm sep}$ is a separable closure of $E$, Fe$(E)$ is the
set of finite extensions of $E$ in $E _{\rm sep}$,
$\mathcal{G}_{E} := \mathcal{G}(E _{\rm sep}/E)$ is the absolute
Galois group of $E$; $N(E _{1}/E)$ denotes the norm group of the
extension $E _{1}/E$, for any $E _{1} \in {\rm Fe}(E)$;

\item For each $p \in \mathbb{P}$, $_{p}{\rm Br}(E) = \{b _{p} \in
{\rm Br}(E)\colon \ pb _{p} = 0\}$ is the maximal subgroup of
Br$(E)$ of period dividing $p$, Br$(E) _{p}$ - the $p$-component
of Br$(E)$, Brd$_{p}(E)$ - the Brauer $p$-dimension of $E$,
abrd$_{p}(E)$ - the absolute Brauer $p$-dimension of $E$; also,
$E(p)$ is the maximal $p$-extension of $E$ (in $E _{\rm sep}$),
and cd$_{p}(\mathcal{G}_{E})$ - the cohomological $p$-dimensions
of $\mathcal{G}_{E}$, in the sense of \cite{S1};

\item For any field extension $E ^{\prime }/E$, $I(E ^{\prime }/E)$
denotes the set of intermediate fields of $E ^{\prime }/E$, and
Br$(E ^{\prime }/E)$ is the relative Brauer group of $E ^{\prime
}/E$;

\item Algebraic structures attached to a field $K$ with a nontrivial
Krull valuation $v$: $O _{v}(K) = \{a \in K\colon \ v(a) \ge 0\}$
- the valuation ring of $(K, v)$; $M _{v}(K) = \{\mu \in K\colon \
v(\mu ) > 0\}$ - the maximal ideal of $O _{v}(K)$; $O _{v}(K)
^{\ast } = \{u \in K\colon \ v(u) = 0\}$ - the multiplicative
group of $O _{v}(K)$; $v(K)$ - the value group of $(K, v)$;
$\overline {v(K)}$ - a divisible hull of $v(K)$; for each $\gamma \in
\overline {v(K)}$, $\gamma \ge 0$, $\nabla _{\gamma }(K)$ denotes the
set $\{\lambda \in K\colon \ v(\lambda - 1) > \gamma \}$;

\item $\widehat K = O _{v}(K)/M _{v}(K)$ is the residue field of $(K,
v)$, and for any $\lambda \in O _{v}(K)$, $\hat \lambda \in
\widehat K$ is the residue class $\lambda + M _{v}(K)$; $(K, v)$
is said to be of mixed characteristic $(0, p)$ if char$(K) = 0$
and char$(\widehat K) = p > 0$;
\label{k999}
\item When $(K, v)$ is a real-valued field, $K _{v}$ stands for the
completion of $K$ with respect to the topology induced by $v$, and
$\bar v$ is the valuation of $K _{v}$ continuously extending $v$;
\label{approx}
\item Given an HDV-field $(K, v)$ of mixed characteristic $(0,
p)$, and a primitive $p$-th root of unity $\varepsilon \in K _{\rm
sep}$, we write $\beta \approx \beta '$, for some $\beta , \beta '
\in K(\varepsilon ) ^{\ast }$, if $v(\beta - \beta ') > p\kappa $,
where $\kappa = v(p)/(p - 1)$; given an element $\pi \in K$ with
$0 < v(\pi ) < p\kappa $, we write $\beta \sim \beta '$ if
$v(\beta - \beta ') > p\kappa - v(\pi ) = v((1 - \varepsilon )
^{p}\pi ^{-1})$.
\end{itemize}

\section{\bf Statement of the main result}

\par
\medskip
Let $(K, v)$ be an HDV-field with char$(\widehat K) = p > 0$. As
shown in \cite{PS}, if char$(K) = 0$ and $[\widehat K\colon
\widehat K ^{p}] = p ^{n}$, for some $n \in \mathbb{N}$, then
$[n/2] \le {\rm abrd}_{p}(K) \le 2n$; abrd$_{p}(K) = \infty $ if
and only if $[\widehat K\colon \widehat K ^{p}] = \infty $ (this
is contained in \cite{PS}, Corollary~2.5 and Lemma~2.6). When
$[\widehat K\colon \widehat K ^{p}] = p ^{n}$ and $n$ is odd, it
has been proved in \cite{BH} that abrd$_{p}(K) \ge 1 + [n/2]$. The
proofs of these results show their validity for Brd$_{p}(K)$ if
$K$ contains a primitive $p$-th root of unity (see Remark
\ref{rema6.2}).
\par
\medskip
The purpose of the present paper, in the first place, is to prove
the inequality Brd$_{p}(K) \ge n$ in general, and thereby, to
obtain the inequality abrd$_{p}(K) \ge n$ in Conjecture
\ref{conj1.1}. Also, its major objective is to give an optimal
infinitude criterion for Brd$_{p}(K)$. Our main result can be
stated as follows:
\par
\medskip
\begin{theo}
\label{theo2.1}
Let $(K, v)$ be an {\rm HDV}-field with {\rm char}$(\widehat K) = p >
0$. Then:
\par
{\rm (a)} {\rm Brd}$_{p}(K)$ is infinite if and only if $\widehat
K/\widehat K ^{p}$ is an infinite extension;
\par
{\rm (b)} There exists $D \in d(K)$ with {\rm exp}$(D) = p$ and
{\rm deg}$(D) = p ^{n}$, provided that $[\widehat K\colon \widehat
K ^{p}] = p ^{n}$, for some $n \in \mathbb{N}$; in particular,
{\rm Brd}$_{p}(K) \ge n$.
\end{theo}
\par
\medskip
Theorem \ref{theo2.1} (b) and the right-to-left implication in
Theorem \ref{theo2.1} (a) are proved in Section 6. For our proof,
we construct in Section 3 an algebra $D \in d(K)$ with exp$(D) =
p$ and deg$(D) = p ^{\mu }$, assuming that $K$ has a TR and Galois
extension $M _{\mu }$ of degree $[M _{\mu }\colon K] = p ^{\mu }
\le [\widehat K\colon \widehat K ^{p}]$ with an abelian Galois
group $\mathcal{G}(M _{\mu }/K)$ of period $p$. By a TR-extension
of $K$, we mean here a finite extension $M/K$ with $\widehat M =
\widehat K$. This agrees, by Lemma \ref{lemm3.2} (b), with the
definition of a TR-extension over any valued field, given before
the statement of Lemma \ref{lemm3.3} (for the case of a discrete
valued field, see the paragraph before the statement of Lemma
\ref{lemm5.2}). The existence of $M _{\mu }$ is a consequence of
the following result, which is of independent interest when $(K,
v)$ is of mixed characteristic $(0, p)$:
\par
\medskip
\begin{lemm}
\label{lemm2.2}
Let $(K, v)$ be an {\rm HDV}-field with {\rm char}$(\widehat K) = p
> 0$ and $\widehat K$ infinite. Then $K$ has {\rm TR} extensions $M
_{\mu }$, $\mu \in \mathbb{N}$, such that $[M _{\mu }\colon K] = p
^{\mu }$, $M _{\mu }/K$ is a Galois extension and the group
$\mathcal{G}(M _{\mu }/K)$ is abelian of period $p$, for each $\mu $.
\end{lemm}
\par
\smallskip
Theorem \ref{theo2.1} and Lemma \ref{lemm2.2} have already been
proved in case char$(K) = p$ (cf. \cite{Ch4}, Lemma~4.2).
Moreover, it follows that, in the setting of the lemma, if
char$(K) = p$, then each finite $p$-group $G$ is isomorphic to
$\mathcal{G}(M _{G}/K)$, for some TR and Galois extension $M _{G}$
of $K$ (see \cite{Ch6}, Lemma~2.3). When char$(K) = 0$ and $(K,
v)$ is an HDV-field of type II, in the sense of Kurihara, this is
not true, for any cyclic $p$-group $G$ of sufficiently large order
\cite{MKu}, 12.2, Theorem~(b).
\par
\smallskip
Lemma \ref{lemm2.2} is proved in Sections 5 and 6. Section 3
contains valuation-theoretic preliminaries used in the sequel. We
also show there how Theorem \ref{theo2.1} can be deduced from
Lemma \ref{lemm2.2} (see Lemma \ref{lemm3.6}). For reasons noted
above, here we focus our attention on the mixed characteristic
case $(0, p)$. For the proof of Lemma \ref{lemm2.2}, we take into
consideration whether or not $v(p) \in pv(K)$ (see Lemmas
\ref{lemm4.8}, \ref{lemm6.1} and Lemma \ref{lemm5.2} (b),
respectively). Section 4 is devoted to the technical preparation
for the proof of Lemma \ref{lemm2.2}. As noted above, in Section
7, we prove Conjecture \ref{conj1.1} for an HDV-field of
characteristic $p$. Open questions concerning Brd$_{p}(K)$ are
also posed in two frequently considered special cases.
\par
\medskip
The basic notation, terminology and conventions kept in this paper
are standard and essentially the same as in \cite{L}, \cite{TW}
and \cite{Ch4}. Missing definitions concerning central simple
algebras can be found in \cite{P}. Throughout, Brauer and value
groups are written additively, Galois groups are viewed as
profinite under the Krull topology, and by a profinite group
homomorphism, we mean a continuous one. For any discrete valued
field $(K, v)$, we suppose that $v(K)$ is chosen to be a subgroup
of the additive group $\mathbb{Q}$ of rational numbers. By an
$n$-dimensional local field, for some $n \in \mathbb{N}$, we mean
a complete $n$-discretely valued field $K _{n}$, in the sense of
\cite{F1} (see also \cite{Zh}), with a quasifinite $n$-th residue
field $K _{0}$.

\medskip
\section{\bf Preliminaries}

\medskip
Let $K$ be a field with a (nontrivial) Krull valuation $v$. We say
that $v$ is Henselian, if it extends uniquely, up-to equivalence,
to a valuation $v _{L}$ on each algebraic extension $L$ of $K$.
This holds, if $K = K _{v}$ and $(K, v)$ is a real-valued field,
i.e. $v(K)$ is isomorphic to an ordered subgroup of the additive
group $\mathbb{R}$ of real numbers (cf. \cite{L}, Ch. XII).
Maximally complete fields are also Henselian, since Henselizations
of valued fields are their immediate extensions (see \cite{E3},
Theorem~15.3.5). The valuation $v$ is Henselian if and only if any
of the following two equivalent conditions holds (cf. \cite{E3},
Sect. 18.1, and \cite{Wa}, Theorem~32.19):
\par
\medskip\noindent
(3.1) (a) Given a polynomial $f(X) \in O _{v}(K) [X]$ and an element
$a \in O _{v}(K)$, such that $2v(f ^{\prime }(a)) < v(f(a))$, where
$f ^{\prime }$ is the formal derivative of $f$, there is a zero $c
\in O _{v}(K)$ of $f$ satisfying the equality $v(c - a) = v(f(a)/f
^{\prime }(a))$;
\par
(b) For each normal extension $\Omega /K$, $v ^{\prime }(\tau (\mu ))
= v ^{\prime }(\mu )$ whenever  $\mu \in \Omega $, $v ^{\prime }$ is
a valuation of $\Omega $ extending $v$, and $\tau $ is a
$K$-automorphism of $\Omega $.
\par
\medskip
When $(K, v)$ is real-valued, it is Henselian if and only if $K$ is
(relatively) separably closed in $K _{v}$ (cf. \cite{E3},
Theorems~15.3.5, 17.1.5). The following lemma allows to extend to the
Henselian case results on complete real-valued fields (e.g., the
Grunwald-Wang theorem, see \cite{LR} and Remark \ref{rema5.3}).
\par
\medskip
\begin{lemm}
\label{lemm3.1} Let $(K, v)$ be a real-valued field, $\bar v$ the
continuous prolongation of $v$ on $K _{v}$, and $(\mathcal{K}, v')$
an intermediate valued field of $(K _{v}, \bar v)/(K, v)$. Suppose
that $(\mathcal{K}, v')$ is Henselian, identify $\mathcal{K} _{\rm
sep}$ with its $\mathcal{K}$-isomorphic copy in $K _{v,{\rm sep}}$,
and let $f$ be the mapping {\rm Fe}$(\mathcal{K}) \to {\rm Fe}(K
_{v})$, by the rule $\Lambda ^{\prime } \to \Lambda ^{\prime }K
_{v}$. Then:
\par
{\rm (a)} $\mathcal{K} _{\rm sep} \cap K _{v} = \mathcal{K}$, and
each $\Lambda \in {\rm Fe}(K _{v})$ contains a primitive element
$\lambda \in \mathcal{K} _{\rm sep}$ over $K _{v}$, such that $[K
_{v}(\lambda )\colon K _{v}] = [\mathcal{K}(\lambda )\colon
\mathcal{K}]$;
\par
{\rm (b)} $\mathcal{K} _{\rm sep}K _{v} = K _{v,{\rm sep}}$ and
$\mathcal{G}_{\mathcal{K}} \cong \mathcal{G}_{K _{v}}$;
\par
{\rm (c)} The correspondence $f$ is bijective and degree-preserving;
moreover, $f$ and the inverse mapping $f ^{-1}: {\rm Fe}(K _{v}) \to {\rm
Fe}(\mathcal{K})$, preserve the Galois property and the isomorphism
class of the corresponding Galois groups;
\par
{\rm (d)} For each $\nu \in \mathbb{N}$ not divisible by {\rm
char}$(K)$, $K _{v} ^{\ast \nu } \cap \mathcal{K} ^{\ast } =
\mathcal{K} ^{\ast \nu }$.
\end{lemm}
\par
\smallskip
\begin{proof}
The conditions on $(K, v)$ and $(\mathcal{K}, v')$ ensure that
$\mathcal{K} _{\rm sep} \cap K _{v} = \mathcal{K}$. The latter
part of Lemma \ref{lemm3.1} (a) can be deduced from Krasner's
lemma (see \cite{L2}, Ch. II, Propositions~3, 4). Lemma
\ref{lemm3.1} (c) follows from Lemma \ref{lemm3.1} (a) and Galois
theory (cf. \cite{L}, Ch. VI, Theorem~1.12), and Lemma
\ref{lemm3.1} (b) - from Lemma \ref{lemm3.1} (a), (c) and the
definition of the Krull topology on $\mathcal{G}_{\mathcal{K}}$
and $\mathcal{G}_{K _{v}}$. Lemma \ref{lemm3.1} (d) is implied by
the density of $\mathcal{K}$ in $K _{v}$, and by the fact that the
set $\nabla _{\gamma }(K _{v}) = \{\alpha \in K _{v}\colon \bar
v(\alpha - 1) > \gamma \}$ is an open subgroup of $K _{v} ^{\ast
\nu }$, provided $\gamma \in \mathbb{R}$ is sufficiently large
(one may put $\gamma = 0$ if char$(\widehat K) \nmid \nu $).
\end{proof}
\par
\medskip\noindent
When $v$ is Henselian, so is $v _{L}$, for any algebraic field
extension $L/K$; in this case, $\widehat L/\widehat K$ is
algebraic as well. We write $v$ instead of $v _{L}$ and view
$v(L)$ as an ordered subgroup of a fixed divisible hull $\overline
{v(K)}$. This is allowed, since $v(K)$ is an ordered subgroup of
$v(L)$, such that $v(L)/v(K)$ is a torsion group; hence, $v(L)$
embeds in $\overline {v(K)}$ as an ordered subgroup. These facts
follow from Ostrowski's theorem (see \cite{E3}, Theorem~17.2.1),
namely, the assertion that if $[L\colon K]$ is finite, then
$[\widehat L\colon \widehat K]e(L/K)$ divides $[L\colon K]$ and
$[L\colon K][\widehat L\colon \widehat K] ^{-1}e(L/K) ^{-1}$ has
no divisor $p \in \mathbb P$, $p \neq {\rm char}(\widehat K)$;
here $e(L/K)$ denotes the ramification index of $L/K$ (the index
$\vert v(L)\colon v(K)\vert $ of $v(K)$ in $v(L)$). We state below
several known criteria that $[L\colon K] = [\widehat L\colon
\widehat K]e(L/K)$:
\par
\medskip
\begin{lemm}
\label{lemm3.2}
Let $(K, v)$ be a Henselian field and $L/K$ a finite extension. Then
$[L\colon K] = [\widehat L\colon \widehat K]e(L/K)$ in the following
cases:
\par
{\rm (a)} If char$(\widehat K) \nmid [L\colon K]$ (apply Ostrowski's
theorem);
\par
{\rm (b)} If $(K, v)$ is HDV and $L/K$ is separable (see \cite{E3}, Sect.
17.4);
\par
{\rm (c)} When $(K, v)$ is maximally complete (cf. \cite{Wa},
Theorem~31.21).
\par\noindent
Under the hypotheses of (c), if {\rm char}$(K) = p > 0$, then $K
^{p}$ is maximally complete (relative to the valuation induced by
$v$) with a residue field $\widehat K ^{p}$ and a value group
$pv(K)$; this ensures that $[K\colon K^{p}]$ is finite if and only if
so are $[\widehat K\colon \widehat K ^{p}]$ and the quotient group
$v(K)/pv(K)$.
\end{lemm}
\par
\medskip\noindent
Assume that $(K, v)$ is a nontrivially valued field. A finite
extension $R$ of $K$ is said to be inertial with respect to $v$, if $R$
has a unique (up-to equivalence) valuation $v _{R}$ extending $v$,
the residue field $\widehat R$ of $(R, v _{R})$ is separable over
$\widehat K$, and $[R\colon K] = [\widehat R\colon \widehat K]$;
$R/K$ is called a TR-extension with respect to $v$, if $v$ has a
unique prolongation $v _{R}$ on $R$, and the index $\vert v
_{R}(R)\colon v(K)\vert $ equals $[R\colon K]$. When $v$ is
Henselian, $R/K$ is TR if and only if $e(R/K) = [R\colon K]$.
Inertial extensions of Henselian fields have useful properties, some
of which are presented by the following lemma (for a proof, see
\cite{TW}, Theorem~A.23):
\par
\medskip
\begin{lemm}
\label{lemm3.3}
Let $(K, v)$ be a Henselian field. Then:
\par
{\rm (a)} An inertial extension $R ^{\prime }/K$ is Galois if and
only if $\widehat R ^{\prime }/\widehat K$ is Galois. When this
holds, $\mathcal{G}(R ^{\prime }/K)$ and $\mathcal{G}(\widehat R
^{\prime }/\widehat K)$ are canonically isomorphic.
\par
{\rm (b)} The compositum $K _{\rm ur}$ of inertial extensions of $K$
in $K _{\rm sep}$ is a Galois extension of $K$ with $\mathcal{G}(K
_{\rm ur}/K) \cong \mathcal{G}_{\widehat K}$.
\par
{\rm (c)} Finite extensions of $K$ in $K _{\rm ur}$ are inertial,
and the natural mapping of $I(K _{\rm ur}/K)$ into $I(\widehat K
_{\rm sep}/\widehat K)$ is bijective.
\end{lemm}
\par
\medskip
It is known (cf. \cite{Sch}, Ch. 2, Sect. 7, and \cite{TW}, Sect.
1.2.2) that if $(K, v)$ is Henselian, then $v$ extends on each $D \in
d(K)$ to a unique valuation $v _{D}$, up-to equivalence. Put $v(D) =
v _{D}(D)$ and denote by $\widehat D$ the residue division ring of
$(D, v _{D})$. Note that $\widehat D$ is a division $\widehat
K$-algebra with $[\widehat D\colon \widehat K] < \infty $, and $v(D)$
is an ordered abelian group including $v(K)$ is an ordered subgroup
of finite index $e(D/K)$. In addition, the following holds, by \cite{TY},
Proposition~2.2:
\par
\medskip
\begin{lemm}
\label{lemm3.4} If $(K, v)$ is an {\rm HDV}-field, then $[D\colon
K] = [\widehat D\colon \widehat K]e(D/K)$, for every $D \in d(K)$.
\end{lemm}
Next we state results on any HDV-field $(K, v)$ that are used in
Section 7 for proving Conjecture \ref{conj1.1} in the case of
char$(K) = p$. They reduce the proof of the upper bound in this
conjecture to considering only the case where $(K, v)$ is a complete
discrete valued field (which allows to apply results of \cite{PS} and
\cite{BH}):
\par
\medskip
\begin{fact}
\label{fact3.5} {\rm (a)} The scalar extension map {\rm Br}$(K)
\to {\rm Br}(K _{v})$ is an injective homomorphism which preserves
Schur indices and exponents (cf. \cite{Cohn}, Theorem~1, and
\cite{Sch}, Ch. 2, Theorem~9); hence, {\rm Brd}$_{p'}(K) \le {\rm
Brd}_{p'}(K _{v})$, for every $p' \in \mathbb P$;
\par
{\rm (b)} The valued field $(K _{v}, \bar v)$ (see page
\pageref{k999}) is maximally complete (cf. \cite{Sch}, Ch. 2,
Theorem~8, or \cite{TW}, Example~3.11); in addition, $(K _{v}, \bar
v)/(K, v)$ is an immediate extension (cf. \cite{E3}, Theorem~9.3.2,
or \cite{L}, Ch. XII, Sect. 5).
\end{fact}
\par
\medskip
Let now $(K, v)$ be an HDV-field with char$(\widehat K) = p$.
Suppose that there exists a Galois extension $M/K$ with
$\mathcal{G}(M/K)$ abelian of period $p$ and order $p ^{\mu }$, for
some $\mu \in \mathbb{N}$. Then, by Galois theory, $M$ equals the
compositum $L _{1} \dots L _{\mu }$ of degree $p$ (Galois) extensions
$L _{j}$ of $K$ in $M$, $j = 1, \dots , \mu $. This enables one to
construct various algebras of degree $p ^{\mu }$ presentable as
tensor products of cyclic $K$-algebras of degree $p$ (concerning
cyclic algebras in general, see, e.g., \cite{P}, Sect. 15). When $M/K$
is a TR-extension and $p ^{\mu } \le [\widehat K\colon \widehat K
^{p}]$, our next lemma provides a criterion for an algebra of this
type to lie in $d(K)$, which is used for proving Theorem
\ref{theo2.1}. Before stating it, note that a finite system $\Theta $
of $m$ elements of a field $E$ with char$(E) = p$ is called
$p$-independent over $E ^{p}$, if $[E ^{p}(\Theta )\colon E ^{p}] = p
^{m}$.
\par
\medskip
\begin{lemm}
\label{lemm3.6}
Let $(K, v)$ be an {\rm HDV}-field with {\rm char}$(\widehat K) = p >
0$, and let $M/K$ be a {\rm TR} and Galois extension with
$\mathcal{G}(M/K)$ abelian of period $p$ and finite order $p ^{\mu }
\le [\widehat K\colon \widehat K ^{p}]$. Fix a presentation $M = L
_{1} \dots L _{\mu }$ as a compositum of degree $p$ extensions of $K$
in $M$, take a generator $\sigma _{j}$ of $\mathcal{G}(L _{j}/K)$,
for each index $j$, and choose elements $a _{j} \in O _{v}(K)$, $j
= 1, \dots , \mu $, so that the system $\hat a _{j} \in \widehat K$, $j
= 1, \dots , \mu $, be $p$-independent over $\widehat K ^{p}$. Then
the tensor product $D _{\mu } = \otimes _{j=1} ^{\mu } \Delta _{j}$
of the cyclic $K$-algebra $\Delta _{j} = (L _{j}/K, \sigma _{j}, a
_{j})$, $j = 1, \dots , \mu $, lies in $d(K)$, where $\otimes =
\otimes _{K}$. Moreover, $v(D _{\mu }) = v(M)$ and $\widehat D _{\mu
}$ is a root field over $\widehat K$ of the binomials $X ^{p} - \hat a
_{j}$, $j = 1, \dots , \mu $, so $[\widehat D _{\mu }\colon \widehat
K] = p ^{\mu }$.
\end{lemm}
\par
\medskip
The proof of Lemma \ref{lemm3.6} is done by induction on $\mu $,
by the method of proving \cite{Ch4}, Lemma~4.2 (b) (which covers
the case of $p = {\rm char}(K)$). For convenience of the reader,
we outline its main steps. In fact, it suffices to prove that $D
_{\mu } \in d(K)$; then the rest of the lemma can be deduced from
Lemma \ref{lemm3.4}, the equality $[D _{\mu }\colon K] = p ^{2\mu
}$, and the existence of $K$-subalgebras $\Theta _{\mu }$ and $W
_{\mu }$ of $D _{\mu }$, such that $\Theta _{\mu } \cong M$ and $W
_{\mu }$ is a root field over $K$ of the binomials $X ^{p} - a
_{j}$, $j = 1, \dots , \mu $. If $\mu = 1$, then $\hat a _{1}
\notin \widehat K ^{p} = \widehat L _{1} ^{p}$, which implies $a
_{1} \notin N(L _{1}/K)$; hence, by \cite{P}, Proposition~15.1~b,
$D _{1} \in d(K)$. When $\mu \ge 2$, it suffices to show that $D
_{\mu } \in d(K)$, under the extra hypothesis that the centralizer
$C = C _{D _{\mu }}(L _{\mu })$ lies in $d(L _{\mu })$. As $C = D
_{\mu -1} \otimes _{K} L _{\mu }$, where $D _{\mu -1} = \otimes
_{j=1} ^{\mu -1} \Delta _{j}$, it is easy to see that $v _{C}(C) =
v(M)$ and $\widehat C$ equals the (commutative) field $\widehat
K(\sqrt[p]{\hat a _{1}}, \dots , \sqrt[p]{\hat a _{\mu -1}})$; in
particular, $\widehat C$ does not possess nontrivial $\widehat
K$-automorphisms. Observing that $D _{\mu } \in s(K)$, consider
the $K$-automorphism $\varphi $ of $C$ which induces the identity
on $D _{\mu -1}$ and the automorphism $\sigma _{\mu }$ of $L _{\mu
}$. It follows from the Skolem-Noether theorem (see \cite{P},
Sect. 12.6) that $\varphi $ is induced by the inner automorphism
of $D _{\mu }$ defined by conjugation by an element $x _{\mu } \in
\Delta _{\mu }$ that induces $\sigma _{\mu }$ on $L _{\mu }$,
satisfies $x _{\mu } ^{p} = a _{\mu }$, and generates $D _{\mu }$
over $C$. Thus, $D _{\mu }$ is a cyclic generalized crossed
product over $C$ as described in \cite{A2}, Ch. XI, Theorems~10,
11. In view of (3.1) (b), it is easily verified that the
composition $v _{C} \circ \varphi $ is a valuation of $C$
extending the prolongation of $v$ on $L _{\mu }$. As $v$ is
Henselian, this means that $v _{C} \circ \varphi = v _{C}$ which
implies $v _{C}(d) = 0$ and $\hat d = \hat d ^{\prime p} \in
\widehat C ^{p}$, provided that $d = \prod _{i=0} ^{p-1} \varphi
^{i}(d ^{\prime })$, for some $d ^{\prime } \in C$ with $v _{C}(d
^{\prime }) = 0$. Since $\hat a _{\mu } \notin \widehat C ^{p}$
(and $v _{C}(d) \neq 0$ if $v _{C}(d ^{\prime }) \neq 0$), one
thereby concludes that $\prod _{i=0} ^{p-1} \varphi ^{i} (\tilde
d) \neq a _{\mu }$, for any $\tilde d \in C$. Hence, by the
equality $x _{\mu } ^{p} = a _{\mu }$ and the hypothesis that $C
\in d(L _{\mu })$, the assertion that $D _{\mu } \in d(K)$ can be
obtained from \cite{A2}, Ch. XI, Theorem~12, so Lemma
\ref{lemm3.6} is proved.
\par
\medskip
Theorem \ref{theo2.1} is implied by Lemmas \ref{lemm3.6} and
\ref{lemm2.2}, so our main goal in the rest of the paper is to prove
Lemma \ref{lemm2.2}. As noted in Section 2, one may consider only the
case of char$(K) = 0$. Our next lemma is used in Section 5 for
proving Lemma \ref{lemm2.2}, under the extra hypothesis that $v(p)
\notin pv(K)$.
\par
\medskip
\begin{lemm}
\label{lemm3.7}
Let $(K, v)/(\Phi , \omega )$ be a valued field extension, such that
the index $\vert v(K)\colon \omega (\Phi )\vert $ of $\omega (\Phi )$
in $v(K)$ is finite, and let $\Psi $ be an extension of $\Phi $ in $K
_{\rm sep}$ of degree $p ^{\mu }$, for some $p \in \mathbb{P}$, $\mu
\in \mathbb{N}$. Suppose that $\Psi $ is {\rm TR} over $\Phi $
relative to $\omega $, and $p \nmid \vert v(K)\colon \omega (\Phi
)\vert $. Then $\Psi K/K$ is {\rm TR} relative to $v$ and $[\Psi
K\colon K] = p ^{\mu }$.
\end{lemm}
\par
\smallskip
\begin{proof}
In view of \cite{E3}, Theorem~15.3.5, and our assumptions, one may
suppose, for the proof, that the value groups of all valuations of
$\Psi K$ extending $\omega $ are ordered subgroups of $\overline
{v(K)}$. Let $v'$ be any valuation of $\Psi K$ extending $v$. By the
Fundamental Inequality (cf. \cite{E3}, Theorem~17.1.5),
\par
\medskip\noindent
(3.2) $\vert v'(\Psi K)\colon v(K)\vert \le [\Psi K\colon K] \le
[\Psi \colon \Phi ] = p ^{\mu }$.
\par
\medskip\noindent
As $\Psi /\Phi $ is TR relative to $\omega $, $\Psi $ has a unique
valuation $\omega '$ extending $\omega $. This shows that $\omega '$
equals the valuation of $\Psi $ induced by $v'$. Note further that
$$p ^{\mu } = \vert \omega '(\Psi )\colon \omega (\Phi )\vert ,$$
$$ \vert \omega '(\Psi )\colon \omega (\Phi )\vert \mid \vert
v'(\Psi K)\colon \omega (\Phi )\vert $$
$${\rm and} \ \vert v'(\Psi K)\colon \omega (\Phi )\vert = \vert
v'(\Psi K)\colon v(K)\vert . \vert v(K)\colon \omega (\Phi )\vert .$$
\par
\medskip\noindent
Since $p \nmid \vert v(K)\colon \omega (\Phi )\vert $ by hypothesis,
it follows that $p ^{\mu } \mid \vert v'(\Psi K)\colon v(K)\vert $,
which
implies the inequalities in (3.2) must be equalities. Hence, $[\Psi
K\colon K] = p ^{\mu }$, and by the Fundamental Inequality, it turns
out that $v'$ is the unique valuation of $\Psi K$ extending $v$ and,
moreover, $\Psi K/K$ is TR relative to $v$, as required.
\end{proof}
\par
\medskip
The next lemma presents well-known properties of binomial
extensions of prime degree, and of cyclotomic extensions. They are
often used without an explicit reference (for a proof of the
lemma, see \cite{L}, Ch. VI, Sects. 3, 9).
\par
\medskip
\begin{lemm}
\label{lemm3.8}
Let $E$ be a field and $p \in \mathbb{P}$. Then:
\par
{\rm (a)} For any $\theta \in E^{\ast }$, the polynomial $X ^{p} -
\theta $ is irreducible over $E$ if and only if it has no root in $E$.
\par
{\rm (b)} If $L/E$ is a finite extension, such that $p \nmid [L\colon
E]$, then $L ^{\ast p} \cap E ^{\ast } = E ^{\ast p}$.
\par
{\rm (c)} If $p \neq {\rm char}(E)$ and $\varepsilon $ is a primitive
$p$-th root of unity in $E _{\rm sep}$, then $E(\varepsilon )/E$ is a
Galois extension with $\mathcal{G}(E(\varepsilon )/E)$ cyclic and
$[E(\varepsilon )\colon E] \mid p - 1$; in particular, $E(\varepsilon
) ^{\ast p} \cap E ^{\ast } = E ^{\ast p}$.
\end{lemm}
\par
\medskip
At the end of this section we recall some known properties of
cyclotomic extensions of valued fields that are used in the sequel.
\par
\medskip
\begin{lemm}
\label{lemm3.9} Let $(K, v)$ be a valued field of mixed
characteristic $(0, p)$ containing a primitive $p$-th root of unity
$\varepsilon $. Then:
\par
{\rm (a)} $v(1 - \varepsilon ) = v(p)/(p - 1)$;
\par
{\rm (b)} $v(-i + \sum _{j=0} ^{i-1} \varepsilon ^{j}) \ge v(1 -
\varepsilon )$, for each $i \in \mathbb{N}$ not divisible by $p$;
\par
{\rm (c)} $v((1 - \varepsilon ) ^{p-1} + p) \ge v((1 - \varepsilon )
^{p}) = pv(p)/(p - 1)$.
\end{lemm}
\par
\smallskip
\begin{proof}
The assumption on char$(\widehat K)$ ensures that $v(p) > 0$,
$\varepsilon \in O _{v}(K)$ and the residue class $\hat \varepsilon $
equals the unit of $\widehat K$. Therefore, $v(1 - \varepsilon ) >
0$, and by the proof of Proposition~4.1.2 (i) of \cite{Co-Th}, $v(p)
= (p - 1)v(1 - \varepsilon )$, as claimed by Lemma \ref{lemm3.9} (a).
Also, the inequality $v(1 - \varepsilon ) > 0$ implies $\hat e _{i} =
i \neq 0$, for each $i \in \mathbb{N}$ not divisible by $p$, where $e
_{i} = \sum _{j=0} ^{i-1} \varepsilon ^{j}$. Lemma \ref{lemm3.9} (b)
follows from the fact that $\mathbb{Z}[\varepsilon ] \subset O
_{v}(K)$ and $\varepsilon - 1$ divides (in the ring
$\mathbb{Z}[\varepsilon ]$) the elements $e _{i} - i = \sum _{j=0}
^{i-1} (\varepsilon ^{j} - 1)$, $i = 1, \dots , p - 1$. Clearly,
Lemma \ref{lemm3.9} (b) shows that $v \big ((p - 1)! - \prod _{i=1}
^{p-1} e _{i}\big ) \ge v(1 - \varepsilon ),$ which implies together
with the equalities
$$\Phi _{p}(1) = \prod _{i=1} ^{p-1} (1 - \varepsilon ^{i}) = p =
(1 - \varepsilon ) ^{p-1}\prod _{i=1} ^{p-1} e _{i},$$ where $\Phi
_{p}(X) = \sum _{j=0} ^{p-1} X ^{j}$ is the $p$-th cyclotomic
polynomial, that
$$v((p - 1)!(1 - \varepsilon ) ^{p-1} - p) \ge v((1 -
\varepsilon ) ^{p}) = pv(p)/(p - 1).$$ \par\noindent As $(p - 1)!
\equiv -1 ({\rm mod} \ p)$ (Wilson's theorem), this proves Lemma
\ref{lemm3.9} (c).
\end{proof}

\section{\bf Normal elements and radical degree $p$ extensions of
HDV-fields}

\par
\medskip
Our goal in this section is to prepare technically the proof of Lemma
\ref{lemm2.2}. In order to achieve it, we need information on the
algebraic properties of $p$-th roots of elements of $\nabla _{0}(F)$,
for a valued field $(F, v)$ of mixed characteristic $(0, p)$. A part
of this information is contained in the following two lemmas.
\par
\medskip
\begin{lemm}
\label{lemm4.1}
Let $(F, v)$ be a valued field of mixed characteristic $(0, p)$, and
let $\alpha \in F$, $\beta \in F ^{\ast }$ be elements,
such that $(1 + \beta ) ^{p} = 1 + \alpha $ and $v(\alpha ) > 0$. Put
$\eta = \alpha - \beta ^{p} - p\beta $ and $\kappa = v(p)/(p - 1)$.
Then $v(\eta ) \ge v(p) + 2v(\beta )$. Moreover,
\par
{\rm (a)} $v(\alpha ) < p\kappa $ if and only if $v(\beta ) < \kappa
$; when this holds, $v(\beta ) = v(\alpha )/p$ and $v(\beta ^{p} -
\alpha ) > v(\alpha ) = v(\beta  ^{p})$.
\par
{\rm (b)} If $v(\alpha ) = p\kappa $, then $v(\beta ) = \kappa $.
\end{lemm}
\par
\smallskip
\begin{proof}
By Newton's binomial formula, one has
$$(1 + \beta ) ^{p} = 1 + \alpha = 1 + \beta ^{p} + \sum
_{i=1} ^{p-1} {p \brack i} \beta ^{i}.$$
Since $v(\alpha ) > 0$ and char$(\widehat F) = p$, this ensures that
$v(\beta ) > 0$. The binomial
\par\vskip0.05truecm\noindent
formula also shows that $\eta = 0$ if $p =
2$, and $\eta = \sum _{i=2} ^{p-1} {p \brack i} \beta ^{i}$ if $p >
2$. Note
\par\vskip0.05truecm\noindent
further that $v({p \brack i}) = v(p)$, for all  $i < p$,
which implies that, in case $p > 2$, the sequence of values $v({p
\brack i}\beta ^{i})$, $i = 1, \dots , p - 1$, strictly increases.
These facts prove
\par\vskip0.05truecm\noindent
that $v(\eta ) \ge v(p) + 2v(\beta )$. The obtained inequality has the
following consequences, which in turn imply statements (a) and (b) of
Lemma \ref{lemm4.1}:
\par
\medskip\noindent
${\rm (i)} \ {\rm If} \ v(\beta ) < \kappa , \ {\rm then} \ v(\alpha
) = v(\beta ^{p}) < pv(\beta ) < p\kappa $
\par\vskip0.08truecm\noindent
${\rm and} \ v(\alpha - \beta ^{p}) = v(p\beta + \eta ) = v(p\beta )
> v(\beta ^{p});$
\par
\vskip0.08truecm\noindent
${\rm (ii)} \ {\rm If} \ v(\beta ) > \kappa , {\rm then} \ v(\beta
^{p}) > v(p\beta ),$ $ v(\alpha - p\beta ) = v(\beta ^{p} + \eta ) >
v(p\beta ),$
\par
\vskip0.09truecm\noindent
${\rm and} \ v(\alpha ) = v(p\beta ) = v(p) + v(\beta ) > p\kappa $;
\par
\medskip\noindent
${\rm (iii)} \ {\rm If} \ v(\beta ) = \kappa , \ {\rm then} \
v(\beta ^{p}) = v(p\beta ) = p\kappa < v(\eta ), \ {\rm whence},
v(\alpha )\ge p\kappa .$
\end{proof}
\par
\medskip
\begin{lemm}
\label{lemm4.2} Let $(F, v)$ be a valued field of mixed
characteristic $(0, p)$, and there exist $\gamma \in F ^{\ast }$ with
$v(\gamma ) = v(p)/(p - 1) : = \kappa $. Assume that $\alpha \in F$
and $\beta \in F ^{\ast }$ satisfy $v(\alpha ) \ge p\kappa $ and $(1 +
\beta ) ^{p} = 1 + \alpha $, put $\delta = \beta /\gamma $, and
denote by $g$ the polynomial $g(X) = \gamma ^{-p}[(1 + \gamma X) ^{p}
- (1 + \alpha )] \in F[X]$. Then:
\par
{\rm (a)} $g$ is monic of degree $p$, $g(\delta ) = 0$ and $g \in O
_{v}(F)[X]$;
\par
{\rm (b)} The reduction $\hat g \in \widehat F[X]$ of $g$ modulo $M
_{v}(F)$ equals $X ^{p} + \hat cX - \hat d$, where $c
= p/\gamma ^{p-1}$, $\hat c \neq 0$ and $d = \gamma ^{-p}\alpha $;
also, $\hat d \neq 0$ if and only if $v(\alpha ) = p\kappa $.
\end{lemm}
\par
\medskip
\begin{proof}
(a): Evidently, $1 + \beta $ is a root of the binomial $X ^{p} - (1 +
\alpha )$, so
\par\vskip0.12truecm\noindent $h(\beta ) = 0$,
where $h(X) = (X + 1) ^{p} - 1 - \alpha = X ^{p} + (\sum _{i=1}
^{p-1} {p \brack i} X ^{p-i}) - \alpha $. Observing also
that $g(X) = \gamma ^{-p}h(\gamma X) = X ^{p} + \sum _{i=1} ^{p-1}
({p \brack i} /\gamma ^{i}).X ^{p-i} - (\alpha /\gamma ^{p}),$ one
obtains that $g(X)$ is monic of degree $p$ and $g(\delta ) = 0$, as
required by
\par\vskip0.11truecm\noindent Lemma \ref{lemm4.2} (a).
Since $v(\alpha ) \ge p\kappa $, $v(\gamma ) = \kappa $, and $p \mid
{p \brack i}$, $i = 1, \dots , p - 1$, it
\par\vskip0.12truecm\noindent is easily verified
that $v(\alpha /\gamma ^{p}) \ge 0$ and $v({p \brack i} /\gamma ^{i})
= (p - i - 1)\kappa \ge 0$, for $i = 1, \dots , p - 1$, proving that
$g(X) \in O _{v}(F)[X]$.
\par\vskip0.1truecm
(b): The preceding calculations show that $v({p \brack p-1}
/\gamma ^{p-1}) = 0$, and in case $p > 2$, they yield $v({p \brack i}
/\gamma ^{i}) > 0$, $i = 1, \dots , p - 2$. Also, by the assumptions
on $v(\gamma )$ and $v(\alpha )$, there exist $c \in O _{v}(F) ^{\ast
}$ and $d \in O _{v}(F)$, such that $p = \gamma ^{p-1}c$ and $\alpha
= \gamma ^{p}d$. These observations show that $\hat g(X) = X ^{p} +
\hat cX - \hat d \in \widehat F[X]$ and $\hat c \neq 0$. They also
prove that $\hat d \neq 0$ if and only if $v(\alpha ) = p\kappa $, as
claimed.
\end{proof}
\par
\medskip
Our approach to the proof of Lemma \ref{lemm2.2} in the case where
$v(p) \in pv(K)$ relies on the following lemma (which is an extended
version of \cite{TY}, Lemma~2.1).
\par
\medskip
\begin{lemm}
\label{lemm4.3}
Let $(K, v)$ be a Henselian field of mixed characteristic $(0, p)$,
$\varepsilon $ be a primitive $p$-th root of unity in $K _{\rm sep}$,
and $\kappa = v(p)/(p - 1)$. Then:
\par
{\rm (a)} The polynomial $g _{\lambda }(X) = (1 - \varepsilon )
^{-p}[((1 - \varepsilon )X + 1) ^{p} - \lambda ]$ lies in $O
_{v}(K(\varepsilon ))[X]$ and has a root in $K(\varepsilon )$, for
each $\lambda \in \nabla _{p\kappa }(K(\varepsilon ))$; in particular,
$\lambda \in K(\varepsilon ) ^{\ast p}$.
\par
{\rm (b)} $\nabla _{\kappa '}(K) \subset K ^{\ast p}$, in case
$\kappa ' \in v(K)$ and $\kappa ' \ge p\kappa $.
\par
{\rm (c)} For any pair $\lambda _{1} \in \nabla _{0}(K)$, $\lambda
_{2} \in K$, such that $v(\lambda _{1} - \lambda _{2}) > p\kappa $,
the elements $\lambda _{2}$ and $\lambda _{2}\lambda _{1} ^{-1}$ lie
in $\nabla _{0}(K)$ and $K ^{\ast p}$, respectively.
\end{lemm}
\par
\smallskip
\begin{proof}
(a): We have $v(1 - \varepsilon ) = \kappa $ and $v(\lambda - 1)
> p\kappa $, whence, Lemma \ref{lemm4.2} applies to $g _{\lambda
}(X)$ and yields $g _{\lambda }(X) \in O _{v}(K(\varepsilon ))[X]$. Denote by
$\widehat K _{\varepsilon }$ the residue field of $(K(\varepsilon ),
v)$. Lemma \ref{lemm4.2}, combined with Lemma \ref{lemm3.9} (c),
shows that the reduction $\hat g _{\lambda }(X) \in \widehat K
_{\varepsilon }[X]$ of $g _{\lambda }(X)$ modulo $M _{v}(K(\varepsilon
))$ equals the binomial $X ^{p} - X$ ($\hat g _{\lambda }(0) = 0$,
since $v((\lambda - 1)/(1 - \varepsilon ) ^{p}) > 0$). This
implies $\hat g _{\lambda }(X)$ has a simple zero in $\widehat K
_{\varepsilon }$, so it follows from (3.1) (a) that $g _{\lambda }(X)$
has a zero in $O _{v}(K(\varepsilon ))$; hence, $\lambda \in
K(\varepsilon ) ^{\ast p}$.
\par\vskip0.05truecm
(b): Lemmas \ref{lemm3.8} (c) and \ref{lemm4.3} (a) imply $\nabla
_{\kappa '}(K) \subset K(\varepsilon ) ^{\ast p} \cap K ^{\ast } = K
^{\ast p}$.
\par\vskip0.05truecm
(c): Clearly, $\nabla _{0}(K)$ contains $\lambda _{2}$ and $\lambda
_{1} ^{-1}$, and $\lambda _{2}\lambda _{1} ^{-1} = 1 + (\lambda
_{2} - \lambda _{1})\lambda _{1} ^{-1}$
\par\vskip0.07truecm\noindent
lies in $\nabla _{p\kappa }(K(\varepsilon ))$, whence, $\lambda
_{2}\lambda _{1} ^{-1} \in (K(\varepsilon ) ^{\ast p} \cap K ^{\ast
}) = K ^{\ast p}$, as claimed.
\end{proof}
\par
\medskip\noindent
{\bf Definition~1.}
An element $\lambda \in \nabla _{0}(K)$, where $(K, v)$ is an
HDV-field of mixed characteristic $(0, p)$, is
called normal over $K$ (or $K$-normal), if $\lambda \notin K ^{\ast
p}$ and $v(\lambda - 1) \ge v(\lambda ^{\prime } - 1)$, for each
element $\lambda ^{\prime }$ of the coset $\lambda K ^{\ast p}$.
\par
\medskip
When $\lambda \notin K ^{\ast p}$, $\lambda K ^{\ast p}$ contains
$K$-normal elements, as Lemma \ref{lemm4.3} (b) and the cyclicity
of $v(K)$ show that the system $v(\lambda ^{\prime } - 1)$,
$\lambda ^{\prime } \in \lambda K ^{\ast p}$, contains a maximal
element $v(\xi - 1)$ (and $\xi $ is $K$-normal). Our next lemma
characterizes $K$-normal elements. Its conclusions follow from
Lemma \ref{lemm3.1} and \cite{Hyo}, Lemma~(2-16), if $K$ contains
a primitive $p$-th root of unity. Stating the lemma, we use the
implication $pv(p)/(p - 1) \in v(K) \Rightarrow v(p)/(p - 1) \in
v(K)$.
\par
\medskip
\begin{lemm}
\label{lemm4.4} Let $(K, v)$ be an {\rm HDV}-field of mixed
characteristic $(0, p)$, and let $\varepsilon $ be a primitive $p$-th
root of unity in $K _{\rm sep}$. Suppose that $\lambda \in \nabla
_{0}(K)$, put $\pi = \lambda - 1$, $\kappa = v(p)/(p - 1)$, and let
$K ^{\prime }$ be an extension of $K$ in $K _{\rm sep}$ obtained by
adjunction  of a $p$-th root $\lambda ^{\prime }$ of $\lambda $. Then
$\lambda $ is $K$-normal if and only if one of the following three
conditions is fulfilled:
\par
{\rm (a)} $v(\pi ) < p\kappa $ and $v(\pi ) \notin pv(K)$; when
this holds, $K ^{\prime }/K$ is {\rm TR};
\par
{\rm (b)} $v(\pi ) < p\kappa $ and $\pi = \pi _{1} ^{p}a$, for
some $\pi _{1} \in K$, $a \in O _{v}(K) ^{\ast }$ with $\hat a
\notin \widehat K ^{\ast p}$; in this case, $\hat a \in \widehat K
^{\prime p}$ and $\widehat K ^{\prime }/\widehat K$ is purely
inseparable of degree $p$;
\par
{\rm (c)} $v(\pi ) = p\kappa $, and for any $\pi _{1} \in K$
with $v(\pi _{1}) = \kappa $, the polynomial $X ^{p} + \hat bX -
\hat d \in \widehat K[X]$ is irreducible over $\widehat K$, $\hat b$
and $\hat d$ being the residue classes of the elements $b = p/\pi
_{1} ^{p-1}$ and $d = \pi /\pi _{1} ^{p}$, respectively; when this
holds, $K ^{\prime }/K$ is inertial and $K(\sqrt[(p-1)]{-b} \ ) =
K(\varepsilon )$.
\end{lemm}
\par
\medskip
\begin{proof}
Put $\pi ^{\prime } = \lambda ^{\prime } - 1$. The conditions of the
lemma show that $v(\pi ) > 0$ and $\lambda ^{\prime } \in \nabla
_{0}(K ^{\prime })$, i.e. $v(\pi ^{\prime }) > 0$. In view of Lemma
\ref{lemm4.3} (b), one may assume, for the proof, that $v(\pi )
\le p\kappa $. Hence, by Lemma \ref{lemm4.1} (a) and (b) (applied to
$(1 + \pi ^{\prime }) ^{p} = 1 + \pi$), $v(\pi ^{\prime }) \le \kappa
$, where equality holds only in case $v(\pi ) = p\kappa $. Our proof
proceeds in three steps.
\par\vskip0.05truecm
Step 1. Let $v(\pi ) < p\kappa $ and $\pi $ violate both
conditions (a) and (b). Then
\par\vskip0.05truecm\noindent $\lambda = 1 + \pi _{0} ^{p}a _{0} ^{p}
+ \pi _{0} ^{\prime }$, for some $a _{0} \in O _{v}(K) ^{\ast }$ and
$\pi _{0}, \pi _{0} ^{\prime } \in K$, such that $v(\pi _{0} ^{p}) =
v(\pi ) < v(\pi _{0} ^{\prime })$. Therefore, applying Lemma
\ref{lemm4.1} to $(1 - \pi _{0}a _{0}) ^{p}$, one obtains that
$v(\lambda (1 - \pi _{0}a _{0}) ^{p} - 1) > v(\pi ) = v(\lambda -
1)$; hence, $\lambda $ is not $K$-normal.
\par\vskip0.05truecm
Step 2. Assume now that $\pi $ satisfies condition (a) or (b) of
Lemma \ref{lemm4.4}. Then, for each $\tilde \lambda \in \nabla
_{0}(K)$ with $v(\tilde \lambda - 1) > v(\pi )$, the element
$\lambda \tilde \lambda - 1$ has value $v(\lambda \tilde \lambda -
1) = v(\pi )$ and satisfies the same condition as $\pi $.
Moreover, under condition (b), $(\lambda \tilde \lambda - 1)/\pi
_{1} ^{p}$ lies in $O _{v}(K) ^{\ast }$ and its residue class
equals $\hat a$. Observing also that $\tilde \lambda ^{-1} \in
\nabla _{0}(K)$ and $v(\tilde \lambda ^{-1} - 1) = v(\tilde
\lambda - 1)$, one concludes that the $K$-normality of $\lambda $
will be proved, if we show that $\lambda \notin K ^{\ast p}$. The
equality $(1 + \pi ^{\prime }) ^{p} = 1 + \pi = \lambda $ and
Lemma \ref{lemm4.1} (a) imply $v(\pi - \pi ^{\prime p}) > v(\pi )
= v(\pi ^{\prime p}) = pv(\pi ^{\prime })$, proving that $v(\pi )
\in pv(K ^{\prime })$. When $v(\pi ) \notin pv(K)$, this means
that $K ^{\prime }/K$ is TR, $[K ^{\prime }\colon K] = p$ and
$\lambda \notin K ^{\ast p}$. Similarly, it follows from Lemma
\ref{lemm4.1} (a) that if $\pi = \pi _{1} ^{p}a$, where $\pi _{1}
\in K$ and $a \in O _{v}(K) ^{\ast }$, then $\pi ^{\prime } = \pi
_{1}a _{1}$, for some $a _{1} \in O _{v}(K ^{\prime }) ^{\ast }$
with $v(a - a _{1} ^{p}) > 0$; hence, $\hat a _{1} ^{p} = \hat a$,
proving that $\hat a \in \widehat K ^{\prime p}$. This shows that
if $\hat a \notin \widehat K ^{p}$, then $[K ^{\prime }\colon K] =
[\widehat K ^{\prime }\colon \widehat K] = p$, $\widehat K
^{\prime }/\widehat K$ is purely inseparable and $\lambda \notin K
^{\ast p}$. Thus our assumptions on $\pi $ guarantee that, in both
cases, $\tilde \lambda ^{-1}\lambda \notin K ^{\ast p}$, for any
$\tilde \lambda \in \nabla _{0}(K)$ with $v(\tilde \lambda - 1) >
v(\pi )$, which implies $\lambda $ is $K$-normal.
\par
Step 3. Suppose that $v(\pi ) = p\kappa $, take $\pi
_{1} \in K$ so that $v(\pi _{1}) = \kappa $, define $b$ and $d$
as in Lemma \ref{lemm4.4} (c), and put $g(X) = \pi _{1} ^{-p}[(1 +
\pi _{1}X) ^{p} - \lambda ]$. It is easily verified that $v(b) =
v(d) = 0$, $g(\pi ^{\prime }/\pi _{1}) = 0$, and $g(X) \in K[X]$ is
monic; also, it follows from Lemma \ref{lemm3.8} (a) that $g(X)$ is
irreducible over $K$ if and only if $\lambda \notin K ^{\ast p}$. At
the same time, Lemma \ref{lemm4.3} (b) implies $\lambda \notin K
^{\ast p}$ if and only if $\lambda $ is $K$-normal. Note further
that, by Lemma \ref{lemm4.2}, $g(X) \in O _{v}(K)[X]$ and its
reduction $\hat g(X)$ modulo $M _{v}(K)$ equals the trinomial $X
^{p} + \hat bX - \hat d \in \widehat K[X]$. In addition, the
equality $v(b) = 0$ shows that $\hat g(X)$ is separable. Using
(3.1)(a), one also proves that $\hat g(X)$ is irreducible over
$\widehat K$ if and only if $\lambda \notin K ^{\ast p}$. It is now
easy to see that $\lambda $ is $K$-normal if and
only if $K ^{\prime }/K$ is inertial with $[K ^{\prime }:K] = p$.
\par
For the rest of the proof of Lemma \ref{lemm4.4} (c), we assume that
$\lambda \notin K ^{\ast p}$, fix a root $\xi \in K _{\rm sep}$ of
the binomial $b(x) = X ^{p-1} + b$, and put $B = K(\xi )$. We first
show that $[B:K] \mid p - 1$. As char$(\widehat K) = p$, $\widehat K$
contains a primitive $(p - 1)$-th root of unity $\hat
\rho $, and since $v$ is Henselian, (3.1) (a), applied to the
binomial $X ^{p-1} - 1$, shows that $\hat \rho $ can be lifted to
such a root $\rho \in K$. Hence, the fact that $[B:K] \mid p - 1$
follows from Galois theory (cf. \cite{L}, Ch. VI, Theorem~6.2).
\par
Finally, we prove that $B = K(\varepsilon )$. It is easily verified
that $\pi ^{\prime }/(\pi _{1}\xi )$ is a root of the monic
polynomial $h(X) = \xi ^{-p}g(\xi X)$. Observing that $v(\xi ) = 0$,
one obtains from the already noted properties of $g(X)$ that $h(X)
\in O _{v}(B)[X]$ and the reduction $\hat h(X) \in \widehat B[X]$ of
$h(X)$ modulo $M _{v}(B))$ is an Artin-Schreier trinomial. Moreover,
it becomes clear that $\widehat h(X) = \hat \xi ^{-p}\hat g(\hat \xi
X)$, which implies in conjunction with Lemma \ref{lemm3.8} (b) (and
the divisibility of $p - 1$ by $[B\colon K]$) that $\hat g(X)$ and $\hat
h(X)$ are irreducible over $\widehat B$. Hence, by Lemma
\ref{lemm3.3} and the Artin-Schreier theorem (cf. \cite{L}, Ch. VI,
Sect. 6), applied to $\hat h(X)$, $K ^{\prime }B/B$ is an inertial
Galois extension of degree $p$. In view of the definition of $K
^{\prime }$, this proves that $\varepsilon \in B$. Let now $\widehat
K _{\varepsilon }$ be the residue field of $(K(\varepsilon ), v)$,
and set $g _{0}(X) = (1 - \varepsilon ) ^{-p}[(1 + (1 - \varepsilon
)X) ^{p} - \lambda ]$. Then $g _{0}(X)$ is monic, and it follows
from Lemmas \ref{lemm4.2} (a), \ref{lemm3.9} (a) that $g _{0}(\pi
^{\prime }/(1 - \varepsilon )) = 0$ and $g _{0}(X) \in O
_{v}(K(\varepsilon ))[X]$. Moreover, Lemmas \ref{lemm3.8} (a) and
\ref{lemm3.9} (c) imply the reduction $\hat g _{0}(X) \in \widehat K
_{\varepsilon }[X]$ is an Artin-Schreier trinomial irreducible over
$\widehat K _{\varepsilon }$. Lemma \ref{lemm4.2}  (b), applied to
$g(X)$ and $g _{0}(X)$, further indicates that if $c = (1 -
\varepsilon )/\pi _{1}$, then $v(c) = 0$ and $\hat c ^{p-1} = -\hat b
\in \widehat K _{\varepsilon }$. Hence, by (3.1) (a), $b(x)$
has a root in $K(\varepsilon )$. As $K$ contains a primitive $(p -
1)$-th root of unity, this means that all roots of $b(X)$ in $K _{\rm
sep}$ in fact lie in $K(\varepsilon )$. It is now obvious that
$B = K(\varepsilon )$, so Lemma \ref{lemm4.4} is proved.
\end{proof}
\par
\medskip
It follows from Lemmas \ref{lemm3.2} (b) and \ref{lemm4.4} that if
$\alpha \in K$ is normal over $K$, then it is normal over any finite
extension of $K$ of prime-to $p$ degree.
\par
\medskip\noindent
{\bf Definition~2.}
In the setting of Lemma \ref{lemm4.4}, an element $\lambda \in
\nabla _{0}(K)$ is called (u)-normal over $K$, where $(u) \in \{(a),
(b), (c)\}$, if it satisfies condition (u).
\par
\medskip
Next we present Albert's characterization \cite{A1}, Ch. IX,
Theorem~6, of Galois extensions of prime degree different from the
characteristic of the ground field. The characterization is based on
Lemma \ref{lemm3.8} (c).
\par
\medskip
\begin{lemm}
\label{lemm4.5}
Assume that $K$ is an arbitrary field, $\varepsilon $ is a primitive $p$-th
root of unity in $K _{\rm sep}$, for some $p \in \mathbb P \setminus
\{{\rm char}(K)\}$, and $\varphi $ a generator of
$\mathcal{G}(K(\varepsilon )/K)$. Fix an integer $s > 0$ satisfying $\varphi
(\varepsilon ) = \varepsilon ^{s}$, and let $\lambda $ be an element
of $K(\varepsilon ) ^{\ast }$. Then the following conditions are
equivalent:
\par
{\rm (a)} $\lambda \notin K(\varepsilon ) ^{\ast p}$ and $\varphi (\lambda
)\lambda ^{-s} \in K(\varepsilon ) ^{\ast p}$;
\par
{\rm (b)} If $L _{\lambda } ^{\prime } = K(\varepsilon
)(\sqrt[p]{\lambda })$, then $L _{\lambda } ^{\prime }$ contains as
a subfield a Galois extension $L _{\lambda }$ of $K$ of degree $p$
(equivalently, the extension $L _{\lambda } ^{\prime }/K$ is Galois with
$\mathcal{G}(L _{\lambda } ^{\prime }/K)$ cyclic and $[L
_{\lambda } ^{\prime }\colon K] = p[K(\varepsilon )\colon K]$).
\end{lemm}
\par
\medskip
Denote by $K(p, 1)$ the compositum of the extensions of $K$ in $K(p)$
of degree $p$, put $K _{\mathcal{G}} = \{\alpha \in K(\varepsilon )
^{\ast }\colon \ \varphi (\alpha )\alpha ^{-s} \in K(\varepsilon )
^{\ast p}\}$, and fix $\ell \in \mathbb{N}$ so that $s\ell \equiv 1
({\rm mod} \ p)$. Obviously, $K _{\mathcal{G}}$ is a subgroup of
$K(\varepsilon ) ^{\ast }$ including $K(\varepsilon ) ^{\ast p}$.
Note also that $K(p, 1)/K$ is a Galois extension with
$\mathcal{G}(K(p, 1)/K)$ abelian of period $p$; this can be deduced
from Galois theory and the normality of maximal subgroups of
nontrivial finite $p$-groups (see \cite{L}, Ch. I, Sect. 6; Ch. VI,
Theorem~1.14). With this notation, Lemma \ref{lemm4.5} can be
supplemented as follows:
\par
\medskip
\begin{lemm}
\label{lemm4.6} {\rm (a)} There is a bijection $\varrho $ of the set
$\Sigma _{p}$ of
finite extensions of $K$ in $K(p,1)$ upon the set $\mathcal{G}_{p}$
of finite subgroups of $K _{\mathcal{G}}/K(\varepsilon ) ^{\ast p}$, such that
\par\noindent
$\varrho (\Lambda ) \cong \mathcal{G}(\Lambda /K) \cong
\mathcal{G}(\Lambda (\varepsilon )/K(\varepsilon ))$, for each
$\Lambda \in \Sigma _{p}$;
\par
{\rm (b)} For each $\lambda \in K(\varepsilon ) ^{\ast }$, the
product $\Omega (\lambda ) = \prod _{j=0} ^{m-1} \varphi
^{j}(\lambda ) ^{\ell (j)}$ lies in $K _{\mathcal{G}}$, where $m =
[K(\varepsilon )\colon K]$ and $\ell (j) = \ell ^{j}$, $j = 0, \dots , m -
1$.
\end{lemm}
\par
\medskip
\begin{proof}
It follows from Lemma \ref{lemm3.8} (c) and Galois theory (cf.
\cite{L}, Ch. VI, Theorem~1.12) that the mapping $\sigma $ of $\Sigma
_{p}$ into the set $\Sigma _{p} ^{\prime }$ of finite extensions of
$K(\varepsilon )$ in $K(p, 1)(\varepsilon )$, by the rule
$\Lambda \to \Lambda (\varepsilon )$, is bijective with
$\mathcal{G}(\Lambda /K) \cong \mathcal{G}(\Lambda (\varepsilon
)/K(\varepsilon ))$, for each $\Lambda \in \Sigma _{p}$. Moreover,
by Kummer theory and Lemma
\ref{lemm4.5}, there is a bijection $\varrho ': \Sigma _{p} ^{\prime
} \to \mathcal{G}_{p}$, such that $\varrho '(\Lambda ^{\prime })
\cong \mathcal{G}(\Lambda ^{\prime }/K(\varepsilon ))$, for each
$\Lambda ^{\prime } \in \Sigma _{p} ^{\prime }$. Therefore, the
composition $\varrho = \varrho ' \circ \sigma $ has the properties
required by Lemma \ref{lemm4.6} (a).
\par We prove Lemma
\ref{lemm4.6} (b). If $\varepsilon \in K$, then the assertion is
obvious, so we assume that $\varepsilon \notin K$, i.e. $m \ge 2$. It
is easily verified that
$$\varphi (\Omega (\lambda )) = \prod _{j=0} ^{m-1} \varphi
^{j+1}(\lambda ) ^{\ell (j)} = \prod _{j=1} ^{m} \varphi ^{j}(\lambda )
^{\ell (j-1)} = \lambda ^{\ell (m-1)}\prod _{j=1} ^{m-1} \varphi
^{j}(\lambda ) ^{\ell (j-1)},$$ \par\vskip0.15truecm\noindent ${\rm and}
\ \Omega (\lambda ) ^{s} = \Omega (\lambda ^{s}) = \prod _{j=0}
^{m-1} \varphi ^{j}(\lambda ) ^{s.\ell (j)}$, for each $\lambda \in
K(\varepsilon )^{\ast }$. Since
\par\vskip0.14truecm\noindent
$s ^{m} \equiv s\ell \equiv 1 ({\rm mod} \ p)$, it follows that $\ell ^{m}
\equiv 1 ({\rm mod} \ p)$, $s \equiv \ell ^{m-1} ({\rm mod} \ p)$,
\par\vskip0.14truecm\noindent
$${\rm and} \ s.\ell (j) \equiv \ell (j - 1) ({\rm mod} \ p), j = 1, \dots,
m - 1,$$ so our calculations prove that $\varphi (\Omega (\lambda
)).\Omega (\lambda ) ^{-s} \in K(\varepsilon ) ^{\ast p}$, as
claimed.
\end{proof}
\par
\medskip
\begin{rema}
\label{rema4.7} Let $(K, v)$ be an HDV-field of mixed
characteristic $(0, p)$, and let $\varepsilon $ be a primitive
$p$-th root of unity in $K _{\rm sep}$. Then:
\par
{\rm (a)} The existence of a {\rm (c)}-normal element over $K$
ensures that $\varepsilon \in K _{\rm ur}$.
\par
{\rm (b)} It can be deduced from Lemma \ref{lemm4.5} that if
$K(\varepsilon )/K$ is TR and $\varepsilon \notin K$ (this holds,
for example, if $v(p)$ generates $v(K)$), then each Galois extension
$L$ of $K$ of degree $p$ is $K$-isomorphic to $L _{\lambda (L)}$,
for some $\lambda (L) \in K _{\mathcal{G}} \cap \nabla
_{0}(K(\varepsilon ))$.
\par
{\rm (c)} When $\langle v(p)\rangle = v(K)$, we have $\langle v(1
- \varepsilon )\rangle = v(K(\varepsilon ))$, which enables one to
obtain from Lemma \ref{lemm4.5}, the preceding observation and
Lemma \ref{lemm4.4} (applied over $K(\varepsilon )$) that a Galois
extension of $K$ of degree $p$ is either inertial or TR (this is a
special case of Miki's theorem, see \cite{MKu}, 12.2). Moreover,
it turns out that degree $p$ extensions of $K _{\rm ur}$ in $K
_{\rm ur}(p)$ are TR (whereas finite extensions of $K _{\rm ur}$
in $K _{\rm ur}(p)$ need not be TR unless $\widehat K$ is perfect,
see Lemmas \ref{lemm5.4} and \ref{lemm3.2} {\rm (b)}).
\end{rema}
\par
\medskip
We conclude this section with the following lemma. As demonstrated in
Section 6, it makes it possible to turn Lemmas \ref{lemm4.3},
\ref{lemm4.4} and \ref{lemm4.6} (a) into the tools we need for the
proof of Lemma \ref{lemm2.2} in the case where $v(p) \in pv(K)$.
\par
\medskip
\begin{lemm}
\label{lemm4.8} Let $(K, v)$ be an {\rm HDV}-field of mixed
characteristic $(0, p)$. Fix a primitive $p$-th root of unity
$\varepsilon \in K _{\rm sep}$, a generator $\varphi $ of
$\mathcal{G}(K(\varepsilon )/K)$, and some $s \in \mathbb{N}$ so that
$\varphi (\varepsilon ) = \varepsilon ^{s}$. Take any $\alpha \in
K(\varepsilon )$ with $v(\alpha ) > v(p)$, and put $\lambda = 1 +
\alpha $. Then $\varphi (\lambda )\lambda ^{-s} \in K(\varepsilon )
^{\ast p}$ in case $v(\varphi (\alpha ) - s\alpha )
> pv(p)/(p - 1)$. \par\vskip0.12truecm\noindent
This holds, if $\alpha = p(1 - \varepsilon )\xi
^{-1}$, where $\xi \in K ^{\ast }$ with $v(\xi ) < v(p)/(p - 1)$.
\end{lemm}
\par
\medskip
\begin{proof}
Put $\kappa = v(p)/(p - 1)$, and {\it use the relation $\approx $
introduced on page \pageref{approx}}. Since for $j \ge 2$,
$v(\alpha ^{j})
> 2v(p) \ge p\kappa $, Newton's binomial formula shows that
\par\vskip0.11truecm\noindent
$\lambda ^{s} \approx 1 + s\alpha $; hence, $\lambda ^{-s} \approx
1 - s\alpha $. Note also that
\par\vskip0.12truecm\noindent $v(\varphi (\alpha )) = v(\alpha )$
because $v$ is Henselian (apply (3.1) (b)). Thus,
$$\varphi (\lambda )\lambda ^{-s} \approx (1 + \varphi (\alpha ))(1
- s\alpha ) \approx (1 + \varphi (\alpha ) - s\alpha ) \approx
1.$$ Hence, $\varphi (\lambda )\lambda ^{-s} \in K(\varepsilon )
^{\ast p}$, by Lemma \ref{lemm4.3} (a).
\par\vskip0.15truecm
Let now $\alpha = p(1 - \varepsilon )\xi ^{-1}$, where $\xi \in K
^{\ast }$ with $0 < v(\xi ) < \kappa $.
\par\vskip0.171truecm\noindent
Then Lemma \ref{lemm3.9} (b) implies the
following, for each $t \in \mathbb{N}$ not divisible by $p$:
$$v(1 - \varepsilon ^{t} - t(1 - \varepsilon )) = v((1 -
\varepsilon ) \sum _{j=0} ^{t-1} (\varepsilon ^{j} - 1)) \ge
2\kappa .$$
Therefore, $v(\alpha ) = v(p) + \kappa - v(\xi ) > v(p)$ and
\par\vskip0.171truecm\noindent
$v(\varphi (\alpha ) - s\alpha ) = v(p[(1 - \varepsilon ^{s}) - s(1
- \varepsilon )]\xi ^{-1}) \ge v(p) + 2\kappa - v(\xi ) > p\kappa .$
\end{proof}

\medskip
\section{\bf Proof of Lemma \ref{lemm2.2} in case char$(K) = 0$ and 
$v(p) \notin pv(K)$}
\par
\medskip
In this section, we consider degree $p$ cyclic extensions related
to Lemma \ref{lemm4.4} (a) and (b), which allows to prove Lemma
\ref{lemm2.2} and Theorem \ref{theo2.1} in the case where char$(K)
= 0$ and $v(p) \notin pv(K)$. Our starting point is the following
lemma.
\par
\medskip
\begin{lemm}
\label{lemm5.1} Let $(K, v)$ be an {\rm HDV}-field of mixed
characteristic $(0, p)$, and let $\varepsilon \in K _{\rm sep}$ be
a primitive $p$-th root of unity, $\varphi $ a generator of
$\mathcal{G}(K(\varepsilon )/K)$, $s$ and $\ell $ positive
integers, such that $\varphi (\varepsilon ) = \varepsilon ^{s}$
and $s\ell \equiv 1 ({\rm mod} \ p)$. Assume that $[K(\varepsilon
)\colon K] = m$, and $\lambda = 1 + (1 - \varepsilon ) ^{p}\pi
^{-1}$, for some $\pi \in K$ with $0 < v(\pi ) < p\kappa $, where
$\kappa = v(p)/(p - 1)$. Denote by $\bar \lambda $ the element
$\Omega (\lambda )$ defined in Lemma \ref{lemm4.6} {\rm (b)}, and
let $L _{\bar \lambda }$ be the extension of $K$ in $K _{\rm sep}$
associated with $\bar \lambda $ in accordance with Lemma
\ref{lemm4.5} {\rm (b)}. Then:
\par
{\rm (a)} If $v(\pi ) \notin pv(K)$, then $\lambda $ and $\bar
\lambda $ are {\rm (a)}-normal over $K(\varepsilon )$; in addition,
$[L _{\bar \lambda }\colon K] = p$, and $L _{\bar \lambda }/K$ is
both Galois and {\rm TR};
\par
{\rm (b)} If $\pi = \pi _{1} ^{p}a$, where $\pi \in K$, $a \in O
_{v}(K) ^{\ast }$ and $\hat a \notin \widehat K ^{p}$, then $\lambda
$ and $\bar \lambda $ are {\rm (b)}-normal over $K(\varepsilon )$;
also, $L _{\bar \lambda }/K$ is Galois, $[L _{\bar \lambda }\colon
K] = p$ and $\widehat L _{\bar \lambda } = \widehat K(\sqrt[p]{\hat
a})$.
\end{lemm}
\par
\medskip
\begin{proof}
Our assumptions and Lemma \ref{lemm3.8} (c) imply $v(\pi ) \in
pv(K)$ if and only if $v(\pi ) \in pv(K(\varepsilon ))$, and
$v(\lambda - 1) \in pv(K(\varepsilon ))$ if and only if $v(\pi ) \in
pv(K)$. They prove that $\widehat K _{\varepsilon } ^{p} \cap
\widehat K = \widehat K ^{p}$, $\widehat K _{\varepsilon }$ being
the residue field of $(K(\varepsilon ), v)$.
\par\vskip0.11truecm\noindent
Putting $e _{n} = \sum _{\nu =0} ^{n-1} \varepsilon ^{\nu }$, for
each $n \in \mathbb{N}$, one obtains from Lemma \ref{lemm3.9} (a),
(b)
\par\vskip0.11truecm\noindent
that $v(n - \varepsilon ^{u}.e _{n}) \ge v(1 - \varepsilon )$, for any
pair $u, n \in \mathbb{N}$ with $p \nmid n$. Since $p \mid n ^{p} -
n$
\par\vskip0.11truecm\noindent
(by Fermat's little theorem), $v(1 - \varepsilon ) = \kappa $, and
$n ^{p} - e _{n} ^{p} = \prod _{u=0} ^{p-1} (n - \varepsilon ^{u}.e
_{n})$,
\par\vskip0.11truecm\noindent
this shows that $v(e _{n} ^{p} - n) \ge v(p),$ which implies the
following:
\par
\vskip0.25truecm\noindent
(5.1) $v((1 - \varepsilon ^{n}) ^{p} - n(1 - \varepsilon ) ^{p}) \ge
v((1 - \varepsilon ) ^{p}) + v(p) > p\kappa .$
\par
\vskip0.22truecm\noindent
Our proof of Lemma \ref{lemm5.1} also relies on the following facts:
\par
\vskip0.22truecm\noindent
(5.2) (a) $v(\bar \lambda - (1 + m(1 - \varepsilon ) ^{p}\pi ^{-1}))
> v((1 - \varepsilon ) ^{p}\pi ^{-1})$;
\par
\medskip
(b) $v(\bar \lambda - 1) = v(m(1 - \varepsilon ) ^{p}\pi ^{-1}) =
p\kappa - v(\pi )$.
\par
\vskip0.22truecm\noindent The equalities in (5.2) (b) follow from
(5.2) (a) (and the equality $v(m) = 0$ implied by Lemma
\ref{lemm3.8} (c)). To prove (5.2) (a) {\it we use the relation
$\sim $ defined on page \pageref{approx} ($\sim $ depends on $\pi
$)}. As $s\ell \equiv 1 ({\rm mod} \ p)$, the relations below,
where $s(j) = s ^{j}$ and $\ell (j) = \ell ^{j}$, include the
content of (5.2) (a) (and forms of (5.1)):
$$\bar \lambda = \prod _{j=0} ^{m-1} [1 + (1 - \varepsilon ^{s(j)})
^{p}\pi ^{-1}] ^{\ell (j)} \sim 1 + \sum _{j=0} ^{m-1} \ell (j)(1 -
\varepsilon ^{s(j)}) ^{p}\pi ^{-1}$$ $$\sim 1 + \sum _{j=0} ^{m-1}
\ell (j)s(j)(1 - \varepsilon ) ^{p}\pi ^{-1} \sim 1 + m(1 -
\varepsilon ) ^{p}\pi ^{-1}.$$ \noindent Statements (5.2) and
observations at the beginning of our proof imply the former parts of
Lemma \ref{lemm5.1} (a) and (b), so we assume further that
either $v(\pi ) \notin pv(K)$ or $\pi = \pi _{1} ^{p}a$, for some
$\pi _{1} \in K$ and $a \in O _{v}(K) ^{\ast }$ with $\hat a \notin
\widehat K ^{p}$. In the former case, $\lambda $ and $\bar \lambda $
are (a)-normal (over $K(\varepsilon )$), and in the latter one, they
are (b)-normal. Let $L _{\bar \lambda } ^{\prime } = K(\varepsilon ,
\bar \lambda ')$, where $\bar \lambda ' \in K _{\rm sep}$ and $\bar
\lambda '^{p} = \bar \lambda $. The normality of $\bar \lambda $
over $K(\varepsilon )$ ensures that $[L _{\bar \lambda } ^{\prime
}\colon K(\varepsilon )] = p$. Using Lemma \ref{lemm4.4}, one
obtains that: if $\bar \lambda $ is (a)-normal, then $L _{\bar
\lambda } ^{\prime }/K(\varepsilon )$ is TR; when $\bar \lambda $ is
(b)-normal, $\widehat L _{\bar \lambda } ^{\prime }/\widehat K
_{\varepsilon }$ is inseparable of degree $p$ with $\hat a \in
\widehat L _{\bar \lambda } ^{\prime p}$. Also, it follows from
Lemmas \ref{lemm4.5}, \ref{lemm4.6} (b) and the $K(\varepsilon
)$-normality of $\bar \lambda $ that $L _{\bar \lambda } ^{\prime }
= L _{\bar \lambda }(\varepsilon )$, and the extension $L _{\bar
\lambda }$ of $K$ in $L _{\bar \lambda } ^{\prime }$ pointed out in
the statement of Lemma \ref{lemm5.1} is Galois with $[L _{\bar
\lambda }\colon K] = p$. As $[L _{\bar \lambda } ^{\prime }\colon L
_{\bar \lambda }] = m$ and $m \mid p - 1$, these observations prove
the following: $L _{\bar \lambda }/K$ is TR if and only if so is $L
_{\bar \lambda } ^{\prime }/K(\varepsilon )$; $\widehat L _{\bar
\lambda }/\widehat K$ is inseparable of degree $p$ if and only if so
is $\widehat L _{\bar \lambda } ^{\prime }/\widehat K _{\varepsilon
}$. Note finally that $[\widehat L _{\bar \lambda } ^{\prime }\colon
\widehat L _{\bar \lambda }] \mid [L _{\bar \lambda } ^{\prime
}\colon L _{\bar \lambda }]$. This implies together with Lemma
\ref{lemm3.8} (b) that if $\bar \lambda $ is (b)-normal, then $\hat a
\in \widehat L _{\bar \lambda } ^{p}$, which completes our proof.
\end{proof}
\par
\medskip
Lemma \ref{lemm3.6} and our next lemma prove Theorem
\ref{theo2.1} in case char$(K) = 0$ and $v(p) \notin pv(K)$. In this
situation, our proof of the lemma relies on the fact (see \cite{FV},
Ch. 2, (3.6), and \cite{E3}, Theorem~15.3.5) that a finite extension
$E ^{\prime }$ of a discrete valued field $(E, w)$ is TR relative to
$w$ if and only if $E ^{\prime }/E$ has a primitive element $\theta $
whose minimal polynomial $f$ over $E$ is Eisenstein at $w$, i.e. $f$
is monic, all of its coefficients but the leading one lie in $M
_{w}(E)$, and the free coefficient of $f$ generates $M _{w}(E)$ as an
ideal of $O _{w}(E)$.
\par
\medskip
\begin{lemm}
\label{lemm5.2}
Let $(K, v)$ be an {\rm HDV}-field of mixed characteristic $(0, p)$.
Suppose that one of the following two conditions is satisfied:
\par
{\rm (a)} $\widehat K$ is an infinite perfect field;
\par
{\rm (b)} $\widehat K$ is imperfect and $v(p) \notin pv(K)$.
\par\noindent
Then there exist {\rm TR} and Galois extensions $M _{\mu }/K$, $\mu
\in \mathbb{N}$, such that $[M _{\mu }\colon K]$ $= p ^{\mu }$ and
$\mathcal{G}(M _{\mu }/K)$ is abelian of period $p$, for each $\mu $.
\end{lemm}
\par
\smallskip
\begin{proof}
We assume, in agreement with conditions (a) and (b), that $\widehat
K$ is infinite. Since the prime subfield, say $\mathbb{F}$, of
$\widehat K$ is finite, this ensures that $\widehat K/\mathbb{F}$ is
an infinite extension, whence, there is a sequence $\tilde b = b
_{\mu } \in O _{v}(K) ^{\ast }$, $\mu \in \mathbb{N}$, such that the
system $\bar b = \hat b _{\mu } \in \widehat K$, $\mu \in
\mathbb{N}$, is linearly independent over $\mathbb{F}$. Denote by
$V$ the $\mathbb{F}$-linear span of the set $\{\hat b _{\mu
}\colon \mu \in \mathbb{N}\}$ and fix a primitive $p$-th root of
unity $\varepsilon \in K _{\rm sep}$, a generator $\varphi $ of
$\mathcal{G}(K(\varepsilon )/K)$, and integers $s$, $\ell $ as in
Lemma \ref{lemm5.1}. Define $K _{\mathcal{G}}$ and $\Omega \colon
K(\varepsilon ) ^{\ast } \to K _{\mathcal{G}}$ as in
Lemma \ref{lemm4.6}, and put $m = [K(\varepsilon )\colon K]$ and
$\lambda _{\mu } = \Omega (1 + (1 - \varepsilon ) ^{p}\pi ^{-1}b
_{\mu })$, $\mu \in \mathbb{N}$, where $\pi \in K$ is fixed so that
$v(\pi ) \notin pv(K)$ and $0 < v(\pi ) \le v(p)$. Take a $p$-th root
$\eta _{\mu } \in K _{\rm sep}$ of $\lambda _{\mu }$, for each $\mu
\in \mathbb{N}$, and consider the fields $L _{\mu } ^{\prime } =
K(\varepsilon , \eta _{\mu })$, $\mu \in \mathbb{N}$. Lemmas
\ref{lemm4.4} and \ref{lemm4.5} show that $[L _{\mu } ^{\prime
}\colon K(\varepsilon )] = p$ and there is a unique Galois extension
$L _{\mu }$ of $K$ in $L _{\mu } ^{\prime }$ of degree $[L _{\mu
}\colon K] = p$. Let $L _{\infty } ^{\prime }$ be the compositum of
the fields $L _{\mu } ^{\prime }$, $\mu \in \mathbb{N}$, and $\Lambda
$ be the subgroup of
$K(\varepsilon ) ^{\ast }$ generated by the set $K(\varepsilon )
^{\ast p} \cup \{\lambda _{\mu }\colon \ \mu \in \mathbb{N}\}$.
Obviously, $\Lambda $ is a subgroup of $K _{\mathcal{G}}$ including
$K(\varepsilon ) ^{\ast p}$. It follows from the assumption on the
sequence $\tilde b$ that, for each $h \in \Lambda \setminus
K(\varepsilon ) ^{\ast p}$, the coset $hK(\varepsilon ) ^{\ast p}$
contains an element of the form $\lambda (h) = 1 + m(1 - \varepsilon
) ^{p}\pi ^{-1}\beta _{h} + \pi (h)$, where $\pi (h) \in
K(\varepsilon )$, $v(\pi (h)) > v(m(1 - \varepsilon ) ^{p}\pi
^{-1})$, $\beta _{h} \in O _{v}(K) ^{\ast }$ and $\hat \beta _{h} \in
V$. Therefore, by the assumptions on $\pi $, $\lambda (h)$ is
(a)-normal over $K(\varepsilon )$ (so Lemma \ref{lemm5.1} (a) applies
to it). This implies $\hat \beta _{h}$ is uniquely determined by $h$
and $\pi $, and does not depend on the choice of $\lambda (h)$ (see
Step 2 of the proof of Lemma \ref{lemm4.4}). More precisely, if $h =
\lambda _{\mu _{1}} ^{k_{1}} \dots \lambda _{\mu _{y}} ^{k _{y}}$,
for some $y \in \mathbb{N}$, and $k _{1}, \dots , k _{y} \in
\mathbb{N}$, with $p \nmid k _{j'}$, for at least one index $j'$,
then $h \notin K(\varepsilon ) ^{\ast p}$, so one may put $\lambda
(h) = h$ and $\beta _{h} = \sum _{j=1} ^{y} k _{j}b _{\mu _{j}}$,
whence, $\hat \beta _{h} = \sum _{j=1} ^{y} k _{j}\hat b _{\mu
_{j}}$. These observations prove that
\par
\vskip0.22truecm\noindent
(5.3) $\{\lambda _{\mu }K(\varepsilon ) ^{\ast p}\colon \mu \in
\mathbb{N}\}$ is a minimal generating set of $\Lambda /K(\varepsilon
) ^{\ast p}$, and there is a unique isomorphism $\rho $ of $\Lambda
/K(\varepsilon ) ^{\ast p}$ upon the additive group of $V$, which
maps the coset $\lambda _{\mu }K(\varepsilon ) ^{\ast p}$ into $\hat
b _{\mu }$, for each $\mu \in \mathbb{N}$.
\par
\vskip0.22truecm\noindent
Statement (5.3), the argument proving it, and Lemmas \ref{lemm4.6}
and \ref{lemm5.1}~(a) imply that the fields $L _{\infty } ^{\prime }$
and $L _{\mu }$, $\mu \in \mathbb{N}$, satisfy the following:
\par
\medskip\noindent
(5.4) (a) $[L _{1} \dots L _{\mu }\colon K] = [L _{1} ^{\prime }
\dots L _{\mu } ^{\prime }\colon K(\varepsilon )] = p ^{\mu }$, for
each $\mu $;
\par
(b) The compositum $L _{\infty }$ of all $L _{\mu }$, $\mu \in
\mathbb{N}$, is an infinite Galois extension of $K$ with $L _{\infty
}(\varepsilon ) = L _{\infty } ^{\prime }$ and $\mathcal{G}(L
_{\infty }/K)$ abelian of period $p$;
\par
(c) Every extension of $K$ in $L _{\infty }$ of degree $p$ is Galois
and TR over $K$.
\par
\medskip
Suppose now that $\widehat K$ is perfect. Then every $R \in {\rm
Fe}(K)$ contains as a subfield an inertial extension $R _{0}$ of $K$
with $\widehat R _{0} = \widehat R$ (cf. \cite{TW},
Proposition~A.17). In view of Lemmas \ref{lemm3.2} and
\ref{lemm3.3} (c), this allows to deduce from (5.4) (b), (c) and
Galois theory that finite extensions of $K$ in $L _{\infty }$ are TR.
Thus the fields $M _{\mu } = L _{1} \dots L _{\mu }$, $\mu \in
\mathbb{N}$, have the properties claimed by Lemma \ref{lemm5.2}.
\par
It remains for us to prove Lemma \ref{lemm5.2} (b). The idea of our
proof has been borrowed from \cite{N}, 2.2.1. Identifying
$\mathbb{Q}$ with the prime subfield of $K$, put $E _{0} =
\mathbb{Q}(t _{0})$, where $t _{0} \in O _{v}(K) ^{\ast }$ is chosen
so that $\hat t _{0} \notin \widehat K ^{p}$ (whence, $\hat t _{0}$
is transcendental over $\mathbb{F}$). Denote by $\omega $ and $v
_{0}$ the valuations induced by $v$ upon $\mathbb{Q}$ and $E _{0}$,
respectively, and fix a system $t _{\mu } \in K _{\rm sep}$, $\mu \in
\mathbb{N}$, such that $t _{\mu } ^{p} = t _{\mu -1}$, for each $\mu
> 0$. It is easy to see that $\mathbb{F}$ equals the residue field of
$(\mathbb{Q}, \omega )$, and the fields $E _{\mu } = \mathbb{Q}(t
_{\mu })$, $\mu \in \mathbb{N}$, are purely transcendental extensions
of $\mathbb{Q}$. Let $v _{\mu }$ be the restricted Gauss valuation of
$E _{\mu }$ extending $\omega $, in the sense of \cite{E3}, for each
$\mu \in \mathbb{N}$. Clearly, for any pair of indices $\nu , \mu $
with $0 < \nu \le \mu $, $E _{\nu - 1}$ is a subfield of $E _{\mu }$
and $v _{\mu }$ is the unique prolongation of $v _{\nu - 1}$ on $E
_{\mu }$. Hence, the union $E _{\infty } = \cup _{\mu =0} ^{\infty }
E _{\mu }$ is a field with a unique valuation $v _{\infty }$
extending $v _{\mu }$, for every $\mu < \infty $. Denote by $\widehat
E _{\mu }$ the residue field of $(E _{\mu }, v _{\mu })$, for each
$\mu \in \mathbb{N} \cup \{0, \infty \}$. The Gaussian property of $v
_{\mu }$, $\mu < \infty $, guarantees that $v _{\mu }(E _{\mu }) =
\omega (\mathbb{Q})$, $v _{\mu}(t _{\mu }) = 0$, $\hat t _{\mu }$ is
a transcendental element over $\mathbb{F}$ and $\widehat E _{\mu } =
\mathbb{F}(\hat t _{\mu })$ (see \cite{E3}, Examples~4.3.2 and
4.3.3). Observing also that $\hat t _{\mu } ^{p} = \hat t _{\mu -1}$,
$\mu \in \mathbb{N}$, $\widehat E _{\infty } = \cup _{\mu =1}
^{\infty } \widehat E _{\mu }$ and $\mathbb{F} ^{p} = \mathbb{F}$,
one concludes that $\widehat E _{\infty }$ is infinite and perfect.
It is therefore clear from
\par\noindent
Lemma \ref{lemm5.2} (a) and Grunwald-Wang's theorem (see Remark
\ref{rema5.3}), that if $(E _{\infty } ^{\prime }, v _{\infty }
^{\prime })$ is a Henselization of $(E _{\infty }, v _{\infty })$
with $E _{\infty } ^{\prime } \subset K _{\rm sep}$, then there exist
TR and Galois extensions $T _{\mu } ^{\prime }/E _{\infty } ^{\prime
}$ and $T _{\mu }/E _{\infty }$, $\mu \in \mathbb{N}$, such that $[T
_{\mu }\colon E _{\infty }] = [T _{\mu } ^{\prime }\colon E _{\infty
} ^{\prime }] = p ^{\mu }$, $T _{\mu } ^{\prime } = T _{\mu }E
_{\infty } ^{\prime }$, $\mathcal{G}(T _{\mu }/E _{\infty })$ is
abelian of period $p$, and $\mathcal{G}(T _{\mu }/E _{\infty }) \cong
\mathcal{G}(T _{\mu } ^{\prime }/E _{\infty } ^{\prime })$, for every
$\mu $. Now fix an arbitrary index $\mu $, choose $\theta \in T _{\mu
}$ so that the minimal polynomial $f(X)$ of $\theta $ over $E
_{\infty }$ be Eisenstein at $v _{\mu }$, and take a sufficiently
large index $k \ge \mu $ such that $f(X) \in E _{k}[X]$ and $f(X)$
splits over $E _{k}(\theta )$. Then the extension $E _{k}(\theta )/E
_{k}$ is both TR and Galois with $\mathcal{G}(E _{k}(\theta )/E _{k})
\cong \mathcal{G}(T _{\mu }/E _{\infty })$, and $f(X)$ is Eisenstein
at $v _{k}$. Let $\psi $ be the isomorphism $E _{k} \to E _{0}$
mapping $t _{k}$ into $t _{0}$. Then $\psi $ extends uniquely to a
degree-preserving isomorphism $\psi ^{\prime }: E _{k}[X] \to E
_{0}[X]$ of polynomial rings, such that $\psi ^{\prime }(X) = X$;
also, $\psi ^{\prime }$ maps $O _{v _{k}}(E _{k})[X]$ into $O _{v
_{0}}(E _{0})[X]$. Note that, for each $g(X) \in E _{k}[X]$,
$\psi ^{\prime }$ induces canonically a ring isomorphism $\psi
^{\prime } _{g}: R _{k} \to R _{0}$ extending $\psi $, where $R _{k}
= E _{k}[X]/(g(X))$ and $R _{0} = E _{0}[X]/(\psi ^{\prime }(g(X)))$.
Clearly, $\psi ^{\prime }_{g}$ maps bijectively the set of roots of
$g(X)$ in $R _{k}$ on the set of roots of $\psi ^{\prime }(g(X))$ in
$R _{0}$. One also sees that $g(X)$ is irreducible over $E _{k}$ if
and only if so is $\psi ^{\prime }(g(X))$ over $E _{0}$. Therefore,
$R _{k}/E _{k}$ is a field extension if and only if so is $R _{0}/E
_{0}$; when this occurs, $[R _{k}\colon E _{k}] = [R _{0}\colon E
_{0}] = {\rm deg}(g)$. Moreover, it follows that $R _{k}/E _{k}$ is
Galois if and only if so is $R _{0}/E _{0}$ (and this holds if and
only if $g(X)$ is irreducible over $E _{k}$ and $R _{k}$ is a root
field of $g(X)$ over $E _{k}$). Suppose now that $R _{k}/E _{k}$ is
Galois. Then, for each $\sigma \in \mathcal{G}(R _{k}/E _{k})$, there
is a unique $\sigma ^{\prime } \in \mathcal{G}(R _{0}/E _{0})$, such
that $\sigma ^{\prime }(\psi ^{\prime } _{g}(r _{k})) = \psi ^{\prime
} _{g}(\sigma (r _{k}))$, for every $r _{k} \in R _{k}$; in addition,
the mapping of $\mathcal{G}(R _{k}/E _{k})$ into $\mathcal{G}(R
_{0}/E _{0})$, by the rule $\sigma \to \sigma ^{\prime }$, is an
isomorphism. Note finally that $v _{k}(e _{k}) = v _{0}(\psi (e
_{k}))$, for every $e _{k} \in E _{k}$, which implies $g(X)$ is
Eisenstein at $v _{k}$ if and only if so is $\psi ^{\prime }(g(X))$
at $v _{0}$; hence, $R _{k}/E _{k}$ is TR relative to $v _{k}$ if and
only if so is $R _{0}/E _{0}$ relative to $v _{0}$. When $g(X) =
f(X)$, these observations show that $\psi $ extends to an isomorphism
of $E _{k}(\theta )$ on the root field $R \in {\rm Fe}(E _{0})$ of
$\psi ^{\prime }(f(X))$ over $E _{0}$, and that $R/E _{0}$ is TR
(relative to $v _{0}$) and Galois with $\mathcal{G}(R/E _{0}) \cong
\mathcal{G}(E _{k}(\theta )/E _{k})$. As $v(p) \notin pv(K)$, one
obtains from Lemma \ref{lemm3.7} and the described properties of $R/E
_{0}$ (regarding $E _{0,{\rm sep}}$ as an $E _{0}$-subalgebra of $K
_{\rm sep}$) that $RK/K$ is TR and Galois, $[RK\colon K] = p ^{\mu }$
and $\mathcal{G}(RK/K) \cong \mathcal{G}(R/E _{0})$ is abelian of
period $p$. Because of the arbitrary choice of the index $\mu $, this
proves Lemma \ref{lemm5.2} (b).
\end{proof}
\par
\smallskip
\begin{rema}
\label{rema5.3} Lemma \ref{lemm3.1} shows that given a field $L$ with
nonequivalent real-valued valuations $w _{1}, \dots , w _{n}$, for
some $n \in \mathbb{N}$, Grunwald-Wang's theorem holds, if applied to
a Henselization of $(L, w _{i})$ (instead of $(L _{w _{i}}, \bar w
_{i})$), for $i = 1, \dots , n$.
\end{rema}
\par
\smallskip
\begin{lemm}
\label{lemm5.4}
Let $(K, v)$ be an {\rm HDV}-field of mixed characteristic $(0, p)$
with $v(p) \in pv(K)$ and $\widehat K \neq \widehat K ^{p}$. Let
$\widetilde \Lambda /\widehat K$ be an inseparable extension of
degree $p$. Then there exists $\Lambda \in I(K(p)/K)$ with $[\Lambda
\colon K] = p$ and $\widehat \Lambda \cong \widetilde \Lambda $ over
$\widehat K$.
\end{lemm}
\par
\smallskip
\begin{proof}
The condition that $v(p) \in pv(K)$ means that there is $\pi _{1} \in
K$ with $v(\pi _{1}) = v(p)/p$, so our conclusion follows
at once from Lemma \ref{lemm5.1} (b).
\end{proof}
\par
\medskip
\begin{prop}
\label{prop5.5} Let $(K, v)$ be an HDV-field of mixed
characteristic $(0, p)$ and with $\widehat K \neq \widehat K
^{p}$. Suppose that $v(p) \in pv(K)$ or $K$ contains a primitive
$p$-th root of unity $\varepsilon $. Then each proper extension
$\widetilde L$ of $\widehat K$ satisfying the inclusion
$\widetilde L ^{p} \subseteq \widehat K$ is $\widehat
K$-isomorphic to $\widehat L$, for some Galois extension $L$ of
$K$, such that $v(L) = v(K)$ and $\mathcal{G}(L/K)$ is an abelian
group of period $p$.
\end{prop}
\par
\medskip
\begin{proof} If $v(p) \in pv(K)$, then our assertion follows from
Lemmas \ref{lemm3.2}, \ref{lemm5.4} and Galois theory; when
$\varepsilon \in K$, it can be deduced from Kummer theory.
\end{proof}
\par
\medskip
\begin{rema}
\label{rema5.6} Let $(K, v)$, $p$ and $\varepsilon \in K _{\rm
sep}$ satisfy the conditions of Lemma \ref{lemm4.4}, and let
$\widehat K \neq \widehat K ^{p}$. Take $c \in O _{v}(K)$ with
$\hat c \notin \widehat K ^{p}$, and suppose that $K$ has a degree
$p$ extension $C$ in $K(p)$, such that $\hat c \in \widehat C
^{p}$. By Lemma \ref{lemm5.4}, an extension of this kind exists if
$v(p) \in pv(K)$ or $\varepsilon \in K$ (this need not hold in
general, see Remark \ref{rema4.7} {\rm (c)}). It is easily
verified that $v(C) = v(K)$, $v(z) \in pv(K)$, for each $z \in
N(C/K)$, and $\hat z \in \widehat C ^{p}$ in case $v(z) = 0$.
Therefore, if $[\widehat K\colon \widehat K ^{p}] \ge p ^{2}$,
$\hat c, \hat b \in \widehat K$ are $p$-independent over $\widehat
K ^{p}$, and $b \in O _{v}(K)$ is a pre-image of $\hat b$, then $b
\notin N(C/K)$ and (by \cite{P}, Proposition~15.1~b) the cyclic
$K$-algebra $V = (C/K, \tau , b)$ of degree $p$ lies in $d(K)$,
$\tau $ being a generator of $\mathcal{G}(C/K)$. Since $v
_{C}(\tau (\alpha ) - \alpha ) > v _{C}(\alpha )$, for any $\alpha
\in C ^{\ast }$, this implies $\widehat V$ contains commuting
$p$-th roots $\hat \eta _{c} = \sqrt[p]{\hat c}$ and $\hat \eta
_{b} = \sqrt[p]{\hat b}$. Hence, by Lemma \ref{lemm3.4}, $v(V) =
v(K)$ and $\widehat V$ equals the field $\widehat K(\hat \eta
_{c}, \hat \eta _{b})$. Also, it follows from Kummer theory that
$V$ is a symbol $K$-algebra, in the sense, e.g. of \cite{PS}, if
and only if $\varepsilon \in K$.
\end{rema}
\par
\medskip\noindent
A detailed and systematic study of algebras $W \in d(K)$, such that
$v(W) = v(K)$ and $\widehat W/\widehat K$ is a purely inseparable
extension would surely be of interest. This, however, goes beyond
the scope of the present paper.
\par
\medskip
\section{\bf Proof of Theorem \ref{theo2.1}}
\par
\medskip
We begin this section with a lemma which completes our preparation
for the proof of Lemma \ref{lemm2.2} in case char$(K) = 0$,
$\widehat K \neq \widehat K ^{p}$ and $v(p) \in pv(K)$. Stating the
lemma, we note that the imposed restriction on $v(p)$ requires the
existence of an element $\pi \in K$ satisfying the conditions
$v(\pi ) \notin pv(K)$ and
$v(p) < pv(p)/(p - 1) - v(p)/p \le v(\pi ) < pv(p)/(p - 1)$.
\par
\medskip
\begin{lemm}
\label{lemm6.1} Let $(K, v)$ be an {\rm HDV}-field of mixed
characteristic $(0, p)$ with $\widehat K$ infinite and $v(p) \in
pv(K)$, and let $\mathbb{F}$ be the prime subfield of $\widehat K$.
Fix an integer $\mu > 0$ and elements $\pi \in K$, $\alpha _{1},
\dots , \alpha _{\mu } \in O _{v}(K) ^{\ast }$, such that $v(\pi )
\notin pv(K)$, $v(p) < v(\pi ) < pv(p)/(p - 1)$, and the system
$\hat \alpha _{1}, \dots , \hat \alpha _{\mu }$ is linearly
independent over $\mathbb{F}$. Put $\lambda _{j} = 1 + \pi \alpha
_{j} ^{p^{\mu }}$, $j = 1, \dots , \mu $, and for any $j$, let $L
_{j} = K(\lambda _{j} ^{\prime })$, where $\lambda _{j} ^{\prime }
\in K _{\rm sep}$ and $\lambda _{j} ^{\prime p} = \lambda _{j}$. Then
the field $M = L _{1} \dots L _{\mu }$ is a {\rm TR}-extension of $K$
of degree $p ^{\mu }$. Moreover, if $K$ contains a primitive $p$-th
root of unity, then $M/K$ is Galois with $\mathcal{G}(M/K)$ abelian
of period $p$.
\end{lemm}
\par
\medskip
\begin{proof}
We first show that one may consider only the special case where $K$
contains a primitive $p$-th root of unity. Let $\varepsilon $ be
such a root in $K _{\rm sep}$. Then $$[M(\varepsilon )\colon K] =
[M(\varepsilon )\colon M][M\colon K] = [M(\varepsilon )\colon
K(\varepsilon )][K(\varepsilon )\colon K].$$ Since, by Galois
theory, $[M(\varepsilon )\colon M] \mid [K(\varepsilon )\colon K]$,
we have $[M(\varepsilon )\colon K(\varepsilon )] \mid [M\colon K].$
In addition, $[K(\varepsilon )\colon K] \mid p - 1$, by Lemma
\ref{lemm3.8}~(c), which implies $p \nmid e(K(\varepsilon )/K)$
(Lemma \ref{lemm3.2}~(b)), proving that $v(\pi ) \notin pv(K(\varepsilon
))$. Moreover, Lemmas \ref{lemm3.2} (b)
and \ref{lemm3.8} imply $\pi $ and $\alpha _{1}, \dots , \alpha _{\mu
}$ satisfy the conditions of Lemma \ref{lemm6.1} with respect to
$(K(\varepsilon ), v)$. Also, it follows from the definition of $M$
that $[M\colon K] \le p ^{\mu }$. As $p \nmid e(M(\varepsilon )/M)$,
these observations prove that if $M(\varepsilon )/K(\varepsilon )$ is
a TR-extension of degree $p ^{\mu }$, then so is $M/K$. This leads to
the desired reduction.
\par
Henceforth, we assume that $\varepsilon \in K$. Then the
concluding assertion of Lemma \ref{lemm6.1} is implied by Kummer
theory and the definition of $M$, so it remains to be seen that
$M/K$ is TR (of degree $p ^{\mu }$). Put $\kappa = v(p)/(p - 1)$
and $\gamma = p\kappa - v(\pi )$. It follows from the conditions
on $\pi $ that $0 < \gamma < \kappa $ and $\gamma \notin pv(K)$.
The rest of our proof relies on the fact that, by Lemma
\ref{lemm4.4} (a), $L _{1}/K$ is TR and $[L _{1}\colon K] = p$,
which means that $M/K$ is TR, provided so is $M/L_{1}$. Using a
standard inductive argument, one may assume for the rest of the
proof that $\mu \ge 2$ and, when $\mu $ is replaced by $\mu - 1$,
the assertion of Lemma \ref{lemm6.1} holds, for any HDV-field $(K
^{\prime }, v ^{\prime })$ of mixed characteristic $(0, p)$ with
$\widehat K ^{\prime }$ infinite and $v ^{\prime }(p) \in pv
^{\prime }(K ^{\prime })$. Then the assertion that $M/L _{1}$ is
TR of degree $p ^{\mu -1}$ can be deduced from the existence of
elements $\pi _{1}$ and $\lambda _{1,j} \in L _{1} ^{\ast }$,
$\alpha _{1,j} \in O _{v}(L _{1}) ^{\ast }$, $j = 2, \dots , \mu
$, such that:
\par\medskip\noindent
(6.1) $\hat \alpha _{1,2}, \dots , \hat \alpha _{1,\mu }$ are
linearly independent over $\mathbb{F}$; $v(\pi _{1}) = p\kappa -
(\gamma /p)$ (whence,
$v(\pi _{1}) \notin pv(L _{1})$); $\lambda _{1,j} = 1 + \pi
_{1}\alpha _{1,j} ^{p ^{\mu -1}}$ and
\par\vskip0.04truecm\noindent $\lambda _{1,j}L _{1} ^{\ast p}
= \lambda _{j}L _{1} ^{\ast p}$, $j = 2, \dots , \mu $.
\par\medskip\noindent
Since the elements $\hat \alpha _{j}\hat \alpha _{1} ^{-1}$,
$j = 1, \dots , \mu $, are linearly independent over $\mathbb{F}$, it
suffices to prove the existence of elements satisfying the conditions
of (6.1) only in the special case where $\alpha _{1} = 1$
(considering $\pi \alpha _{1} ^{p ^{\mu }}$ and $\alpha _{2}\alpha
_{1} ^{-1}, \dots , \alpha _{\mu }\alpha _{1} ^{-1}$
\par\vskip0.04truecm\noindent
instead of $\pi $ and $\alpha _{2}, \dots , \alpha _{\mu }$,
respectively). Putting $\eta _{1} = \lambda _{1} ^{\prime } - 1$, we
show that, in this case, $\pi _{1}$ and $\alpha _{1,j}, \lambda
_{1,j}$, $j = 2, \dots , \mu $, can be chosen as follows:
\par
\medskip\noindent
(6.2) $\pi _{1} = -p\eta _{1}$,  $\alpha _{1,j} = \alpha _{j} -
\alpha _{j} ^{p}$, and $\lambda _{1,j} = 1 - p\eta _{j}$, where $\eta
_{j} = \eta _{1}\alpha _{1,j} ^{p ^{\mu -1}}$.
\par
\medskip\noindent
In the rest of the proof, {\it we use the relation $\approx $
introduced on page \pageref{approx}}. As \par\noindent $(1 + \eta
_{1}) ^{p} = 1 + \pi $ (and $p \ge 2$), Lemma \ref{lemm4.1} (a)
shows that
$$v(\eta _{1}) = v(\pi )/p > v(p)/p, \ {\rm so} \ v(p\eta _{1} ^{2})
> (p + 2)v(p)/p \ge p\kappa ;$$
hence; by the former conclusion of Lemma \ref{lemm4.1}, $\pi \approx
\eta _{1} ^{p} + p\eta _{1}$. At the same time, the equality $\lambda _{j} = 1 +
\pi \alpha _{j} ^{p ^{\mu }}$ implies $\lambda _{j} ^{-1} \approx 1 -
\pi \alpha _{j} ^{p ^{\mu }}$. Likewise, from $\lambda _{1,j} = 1 -
p\eta _{j}$, one obtains that $\lambda _{1,j} ^{-1} \approx 1 + p\eta
_{j}$. Let $\Omega _{j} = 1 + \eta _{1}\alpha _{j} ^{p ^{\mu -1}}$.
Then
$$\Omega _{j} ^{p} \approx 1 + \eta _{1} ^{p}\alpha _{j} ^{p ^{\mu
}} + p\eta _{1}\alpha _{j} ^{p ^{\mu -1}} = [1 + (\eta _{1} ^{p} +
p\eta _{1})\alpha _{j} ^{p ^{\mu }}] + [p\eta _{1}(\alpha _{j} ^{p
^{\mu -1}} - \alpha _{j} ^{p ^{\mu }})]$$ $$\approx \lambda _{j} +
[p\eta _{1}(\alpha _{j} - \alpha _{j} ^{p}) ^{p ^{\mu -1}}] =
\lambda _{j} + p\eta _{j} = \lambda _{j}(1 + p\eta _{j}\lambda _{j}
^{-1}) \approx \lambda _{j}(1 + p\eta _{j}) \approx \lambda
_{j}\lambda _{1,j} ^{-1}.$$
Hence, by Lemma \ref{lemm4.3} (c), $\lambda _{j}\lambda _{1,j} ^{-1}
\in L _{1} ^{\ast p}$.
\par\vskip0.14truecm
We are now in a position to prove Lemma \ref{lemm6.1}. As already
shown,
$$v(p) < v(\pi _{1}) = v(p\eta _{1}) = v(p) + v(\eta _{1}) = p\kappa
- (\gamma /p)$$
and $pv _{1}(L) = v(K)$, which implies $v(\pi _{1}) \notin pv(L
_{1})$. Observing that $\alpha _{1} = 1$,
\par\vskip0.11truecm\noindent
the field $\mathbb{F}$ equals the set $\{\hat y \in \widehat K\colon
\hat y ^{p} = \hat y\}$, ${\rm and} \ \alpha _{1,j} = \alpha _{j} -
\alpha _{j} ^{p}$, $j = 2, \dots , \mu $,
\par\vskip0.08truecm\noindent
are elements of $O _{v}(L _{1}) ^{\ast }$, such that $\hat \alpha
_{1}, \dots , \hat \alpha _{\mu }$ are linearly independent (over
$\mathbb{F}$), one concludes that $\hat \alpha _{1,2}, \dots , \hat
\alpha _{1,\mu }$ are linearly independent as well. Thus the field
$M$ and the elements $\pi _{1} = -p\eta _{1}$, and $\alpha _{1,j},
\lambda _{1,j}$, $j = 2, \dots \mu $, defined in (6.2) satisfy the
conditions of Lemma \ref{lemm6.1} (over $L _{1}$), and by the
inductive hypothesis, $M/L _{1}$ is a TR-extension of degree $p ^{\mu
-1}$, so Lemma \ref{lemm6.1} is proved.
\end{proof}
\par
\medskip
We can now take the final step towards the proof of Lemma
\ref{lemm2.2} (and Theorem \ref{theo2.1}) in general. In view of
Lemma \ref{lemm3.6} and \cite{Ch4}, Lemma~4.2, one may consider only
the case of mixed characteristic $(0, p)$. We also assume that $v(p)
\in pv(K)$ and $\widehat K \neq \widehat K ^{p}$, which is allowed by
Lemma \ref{lemm5.2}. As $v(K)$ is cyclic, the condition on $v(p)$
ensures that there is $\xi \in K$ with $0 < v(\xi ) \le v(p)/p$ and
$v(\xi ) \notin pv(K)$. Since $\widehat K$ is infinite, there are
$\alpha _{\nu } \in O _{v}(K) ^{\ast }$, $\nu \in \mathbb{N}$, such
that the system $\hat \alpha _{\nu } \in \widehat K$, $\nu \in
\mathbb{N}$, is linearly independent over the prime subfield of
$\widehat K$. Take a primitive $p$-th root of unity $\varepsilon \in
K _{\rm sep}$, a generator $\varphi $ of $\mathcal{G}(K(\varepsilon
)/K)$, and $s \in \mathbb{N}$ so that
$\varphi (\varepsilon ) =
\varepsilon ^{s}$. Fix any $\mu \in \mathbb{N}$, put $\lambda _{j} =
1 + p(1 - \varepsilon )\xi ^{-1}\alpha _{j} ^{p ^{\mu }}$, for $j =
1, \dots , \mu $, and denote by $M _{\mu } ^{\prime }$ the
extension of $K(\varepsilon )$ generated by the set $\{\lambda
_{j} ^{\prime }\colon j = 1, \dots , \mu \}$, where $\lambda _{j}
^{\prime } \in K _{\rm sep}$ and $\lambda _{j} ^{\prime p} =
\lambda _{j}$, for any index $j$. It follows from Lemma
\ref{lemm6.1} that $M _{\mu } ^{\prime }/K(\varepsilon )$ is TR and
Galois of degree
$p ^{\mu }$ with $\mathcal{G}(M _{\mu } ^{\prime }/K(\varepsilon
))$ abelian of period $p$. Furthermore, Lemma \ref{lemm4.8} and
the conditions on $\xi $ and $\alpha _{1}, \dots , \alpha _{\mu }$
show that
\par\vskip0.04truecm\noindent $\varphi
(\lambda _{j})\lambda _{j} ^{-s} \in K(\varepsilon ) ^{\ast p}$, $j =
1, \dots , \mu $. Therefore, Lemmas \ref{lemm4.5} and
\ref{lemm4.6}~(a) yield
\par\vskip0.04truecm\noindent
 $M _{\mu } ^{\prime } = M _{\mu }(\varepsilon
)$, for some Galois extension $M _{\mu }$ of $K$ in $K(p)$, such that
\par\vskip0.04truecm\noindent
$\mathcal{G}(M _{\mu }/K) \cong \mathcal{G}(M _{\mu }
^{\prime }/K(\varepsilon ))$; hence, $[M _{\mu }\colon K] = p
^{\mu }$ and $[M _{\mu } ^{\prime }\colon M _{\mu }] = [K(\varepsilon
)\colon K]$. As $p \nmid [K(\varepsilon )\colon K]$ and $M _{\mu }
^{\prime }/K(\varepsilon )$ is TR, it is now easy to see that $M
_{\mu }/K$ is also TR. Because of the arbitrary choice of $\mu
$, this proves Lemma \ref{lemm2.2}, Theorem \ref{theo2.1} (b) and
the right-to-left implication in Theorem \ref{theo2.1} (a).
Finally, by Fact \ref{fact3.5}, the converse implication follows
from \cite{PS}, Corollary~2.5, so Theorem \ref{theo2.1} is proved.
\par
\medskip
\begin{rema}
\label{rema6.2} It should be pointed out that in case $(K, v)$ is
an HDV-field containing a primitive $p$-th root unity $\varepsilon $,
the right-to-left-implication in Theorem \ref{theo2.1} (a) becomes
obvious as a result of the proof of the lower bound for abrd$_{p}(K)$
in \cite{PS}, Lemma~2.6. The conditions of the cited lemma do not
require that $\varepsilon \in K$. However, the assumption that
$\varepsilon \in K$ is necessary to define over $K$ tensor products
of symbol algebras like those used in the proof of \cite{PS},
Lemma~2.6. This allows to show easily that if $\varepsilon \in K$,
then the lower bound in the cited lemma is also such a bound for
Brd$_{p}(K)$, which proves the right-to-left implication in Theorem
\ref{theo2.1} (a).
\end{rema}
\par
\medskip
To end the present Section, we note that Theorem~2 of \cite{PS}
and the conclusion of Theorem \ref{theo2.1} {\rm (b)} in case
char$(K) = 0$ leave open the question of whether abrd$_{p}(E) >
2{\rm Brd}_{p}(E) + 1$, for any field $E$ with a primitive $p$-th
root of unity and Brd$_{p}(E) < \infty $. Moreover, it seems to be
unknown whether abrd$_{p}(E) = \infty $.
\par
\medskip
\section{\bf Open problems and further results}
\par
\medskip
We begin this section with a proof of Conjecture \ref{conj1.1} in
case char$(K) = p$.
\par
\medskip
\begin{prop}
\label{prop7.1} If $(K, v)$ is an {\rm HDV}-field with char$(K) =
p
> 0$, then:
\par
{\rm (a)} Brd$_{p}(K) = \infty $ if $[\widehat K\colon \widehat K
^{p}] = \infty $; when $(K, v)$ is complete, the equality $[\widehat
K\colon \widehat K ^{p}] = \infty $ holds if and only if $[K\colon K
^{p}] = \infty $;
\par
{\rm (b)} $n \le {\rm Brd}_{p}(K) \le n + 1$, provided that $n <
\infty $ and $[\widehat K\colon \widehat K ^{p}] = p ^{n}$;
\par
{\rm (c)} If $(K, v)$ is complete, $[\widehat K\colon \widehat K
^{p}] = p ^{n}$ and $K ^{\prime }/K$ is a finite field extension,
then $[K ^{\prime }\colon K ^{\prime p}] = p ^{n+1}$.
\end{prop}
\par
\medskip\noindent
\begin{proof}
The former part of Proposition \ref{prop7.1} (a) and the lower
bound on Brd$_{p}(K)$ in Proposition \ref{prop7.1} (b) are implied
by \cite{Ch4}, Lemma~4.2 (b). Proposition \ref{prop7.1} (c) and
the latter part of Proposition \ref{prop7.1} (a) follow from Fact
3.5 (b), Lemma \ref{lemm3.2}, and the equality $[L\colon
L ^{p}] = [K\colon K ^{p}]$, for every finite extension $L/K$ (cf.
\cite{BH}, Lemma~2.12). It remains to prove the upper bound in
Proposition \ref{prop7.1} (b). Let $\overline K$ be an algebraic
closure of $K$. In view of \cite{Ch2}, Lemma~4.1, it suffices to
show that, for any finite extension $K ^{\prime }$ of $K$ in
$\overline K$, we have deg$(D ^{\prime }) \mid p ^{n+1}$
whenever $D ^{\prime } \in d(K ^{\prime })$ and exp$(D ^{\prime }) =
p$. In addition, Fact \ref{fact3.5} (a) allows us to consider only
the case of $K = K _{v}$. Let $K _{1} ^{\prime } = \{\lambda \in
\overline K\colon \lambda ^{p} \in K ^{\prime }\}$. Then $K _{1}
^{\prime } \in I(\overline K/K ^{\prime })$, $K _{1} ^{\prime p} =
K ^{\prime }$, and by Proposition \ref{prop7.1} (c), $[K _{1}
^{\prime }\colon K ^{\prime }] = p ^{n+1}$. Since, by Albert's
theorem, $_{p}{\rm Br}(K ^{\prime })$ is a subgroup of Br$(K _{1}
^{\prime }/K ^{\prime })$ (cf. \cite{A2}, Ch. VII, Theorem~28),
this yields deg$(D ^{\prime }) \mid p ^{n+1}$ (see \cite{P}, Sect.
13.4), so Proposition \ref{prop7.1} is proved.
\end{proof}
\par
\medskip
Our next result proves Conjecture \ref{conj1.1} in the special case
where $\widehat K$ is an $n$-dimensional local field of characteristic
$p$ with a finite $n$-th residue field.
\par
\medskip
\begin{prop}
\label{prop7.2}
Assume that $(K, v)$ is an {\rm HDV}-field, such that $\widehat K$ is
an $n$-dimensional local field with {\rm char}$(\widehat K) = p$.
Then {\rm Brd}$_{p}(K) \ge n$. Moreover, if the $n$-th residue field
$\widehat K _{0}$ of $\widehat K$ is finite, then {\rm abrd}$_{p}(K)
\le n + 1$.
\end{prop}
\par
\smallskip
\begin{proof}
As $[\widehat K\colon \widehat K ^{p}] = p ^{n}$, Theorem
\ref{theo2.1} (b) yields Brd$_{p}(K) \ge n$, so it suffices to
prove that if $\widehat K _{0}$ is finite, then abrd$_{p}(K) \le n
+ 1$. In view of Proposition \ref{prop7.1} (b) and Fact
3.5 (a), one may consider only the case of char$(K) = 0$
and $K = K _{v}$. Then $K$ is an $(n + 1)$-dimensional local field
with last residue field $\widehat K _{0}$, whence, by \cite{Ch6},
Proposition~4.4, abrd$_{p}(K) \le n + 1$, as required.
\end{proof}
\par
\medskip
It would be of interest to know whether an HDV-field $(K, v)$ with
$\widehat K _{\rm sep} = \widehat K$ and $[\widehat K\colon
\widehat K ^{p}] = p ^{n}$, for some $n \in \mathbb{N}$, satisfies
Brd$_{p}(K) = n$ (see page \pageref{stconj}). This is the same as to
find whether Brd$_{p}(K) = n$, provided that $p \nmid [\widehat
K^{\prime }\colon \widehat K]$ when $\widehat K ^{\prime }$ runs across
Fe$(\widehat K)$ (cf. \cite{P}, Sects. 13.4 and 14.4). The condition
on $\widehat K$ means that cd$_{p}(\mathcal{G}_{\widehat K}) = 0$. If
$\widehat K _{\rm sep} \neq \widehat K$, then it is possible that
Brd$_{p}(K) \ge n + 1$; such is the case where $\widehat K/\mathbb{F}
_{p}$ is a finitely-generated extension of transcendence  degree $n$
(see the proof of \cite{Ch4}, Proposition~6.3, or \cite{BH},
Theorem~5.2). The same inequality for Brd$_{p}(K)$ is obtained by the
method of proving \cite{Ch4}, Proposition~6.3, when char$(\widehat K)
= p$ and $\widehat K$ is a finitely-generated extension of
transcendence degree $n > 0$ over a perfect field $\widehat K _{0}$
with cd$_{p}(\mathcal{G}_{\widehat K _{0}}) \neq 0$ (see \cite{S1},
Ch. I, 3.3). Since cd$_{p}(\mathcal{G}_{\widehat K _{0}}) \le 1$ (cf.
\cite{S1}, Ch. II, 2.2), Theorem \ref{theo2.1} (b) and the preceding
observations attract interest in the following special case of
Conjecture \ref{conj1.1}:
\par
\smallskip
\begin{conj}
\label{conj7.3} If $(K, v)$ is an HDV-field with char$(\widehat K)
= p > 0$ and $\widehat K$ is a finitely-generated extension of
transcendence degree $n > 0$ over its maximal perfect subfield
$\widehat K _{0}$, then {\rm Brd}$_{p}(K) = n + {\rm
cd}_{p}(\mathcal{G}_{\widehat K _{0}})$.
\end{conj}
\par
\medskip
Theorem \ref{theo2.1} (b) and the upper bounds in \cite{PS},
Theorem~2, \cite{BH}, Corollary~4.7 and Theorem~4.16, and
Proposition \ref{prop7.1} (b) of the present paper prove
Conjecture \ref{conj1.1}, for $n = 1, 2, 3$. Note also that
Conjecture \ref{conj7.3} holds, for $n = 1, 2$. In view of the
remarks preceding the statement of Conjecture \ref{conj7.3}, this
can be obtained by using Theorem \ref{theo2.1} (b), \cite{BH},
Theorem~4.16, and Case IV of the proof of \cite{BH}, Theorem~5.3.
As to Conjecture \ref{conj7.3}, it need not be true if $(K, v)$ is
merely HDV with char$(\widehat K) = p$ and $[\widehat K\colon
\widehat K ^{p}] < \infty $. One may take as a counter-example the
iterated formal power series field $K = \widehat K _{0}((X _{1}))
\dots ((X _{n}))((Y))$  in a system of variables $X _{1}, \dots ,
X _{n}, Y$ over a quasifinite field $\widehat K _{0}$ with
char$(\widehat K _{0}) = p$. Then Brd$_{p}(K) = n$, by \cite{Ch6},
Proposition~3.5 (implied by \cite{Ch2}, Lemma~4.3~(b), \cite{Ch4},
Lemma~4.2 and \cite{AJ}, Theorem~3.3), whereas the formula in
Conjecture \ref{conj7.3} requires Brd$_{p}(K) = n + 1$ (the standard
discrete valuation on $K$ is Henselian with $\widehat K = \widehat K
_{0}((X _{1})) \dots ((X _{n}))$, whence, $[\widehat K\colon \widehat
K ^{p}] = p ^{n}$ and cd$_{p}(\widehat K _{0}) = 1$). This example as
well as Proposition \ref{prop7.2} draw one's attention to the
following problem:
\par
\medskip
\begin{prob}
\label{prob7.4} Let $(K, v)$ be an {\rm HDV}-field with
char$(\widehat K) = p > 0$. Suppose that $\widehat K$ is an
$n$-dimensional local field, for some $n \in \mathbb{N}$, with an
$n$-th residue field $\widehat K _{0}$. Find whether {\rm
Brd}$_{p}(K) = n$.
\end{prob}
\par
\medskip
The conditions of Problem \ref{prob7.4} show that $K _{v}$ is an
$(n + 1)$-dimensional local field with last residue field
$\widehat K _{0}$ (and $\widehat K$ is isomorphic to an iterated
formal power series field in $n$ variables over the quasifinite
field $\widehat K _{0}$, see \cite{F1}, 2.5.2). Therefore, in case
char$(K) = p$, Fact \ref{fact3.5} (a) and \cite{Ch6},
Proposition~3.5, give an affirmative answer to Problem
\ref{prob7.4}. When $n = 1$, such an answer is contained in the
following result of \cite{Ch7}, obtained as a final step towards a
full characterization of stable HDV-fields by properties of their
residue fields:
\par
\medskip
\begin{prop}
\label{prop7.5} Let $(K, v)$ be an {\rm HDV}-field with {\rm
char}$(\widehat K) = p > 0$. Then {\rm Brd}$_{p}(K) \le 1$ if and
only if the following condition is fulfilled:
\par
$[\widehat K\colon \widehat K ^{p}] \le p$, and in case
{\rm Brd}$_{p}(\widehat K) \neq 0$, every degree $p$ extension of
$\widehat K$ in $\widehat K(p)$ is embeddable as a $\widehat
K$-subalgebra in each $D _{p} \in d(\widehat K)$ of degree $p$.
\par\noindent
The equality {\rm Brd}$_{p}(K) = 0$ holds if and only if $\widehat K$
is perfect and $\widehat K(p) = \widehat K$.
\end{prop}
\par
\medskip
\begin{rema}
\label{rema7.6} The inequalities $n \le {\rm Brd}_{p}(K) \le n +
1$ hold, for any HDV-field $(K, v)$, such that $\widehat K$ is an
$n$-dimensional local field with a finite $n$-th residue field and
with char$(\widehat K _{1}) = p$, $\widehat K _{1}$ being the $(n
- 1)$-th residue field of $\widehat K$. Proposition \ref{prop7.2}
reduce the proofs to the case of char$(\widehat K) = 0$ (and $n
\ge 3$, in view of Proposition \ref{prop7.5}). Then the stated
inequalities are contained in \cite{Ch6}, Proposition~4.4.
\end{rema}
\par
\medskip
Note finally that the interest in the question of whether Brd$_{p}(K)
= n$, if $(K, v)$ is an HDV-field, char$(\widehat K) = p > 0$,
$\widehat K _{\rm sep} = \widehat K$ and $[\widehat K\colon \widehat
K ^{p}] = p ^{n}$, for some $n \in \mathbb{N}$, is motivated not only
by Theorem \ref{theo2.1} (b) and \cite{BH}, Theorem~4.16, but also by
the following well-known conjecture (see, e.g., \cite{ABGV}, Sect. 4):
\par
\medskip
\begin{conj}
\label{conj7.7} Assume that $F$ is a field of type $C _{\nu }$,
i.e. each homogeneous polynomial $f(X _{1}, \dots , X _{m}) \in
F[X _{1}, \dots , X _{m}]$ of degree $d$ with $0 < d ^{\nu } < m$,
has a nontrivial zero over $F$. Then abrd$_{p}(F) < \nu $.
\end{conj}
\label{stconj}
\par
\medskip\noindent
To show how Conjecture \ref{conj7.7} is related to the noted
question, fix an HDV-field $(E, \omega )$ so that char$(\widehat
E) = p > 0$, $\widehat E$ be algebraically closed, and when
char$(E) = p$, $E = E _{\omega }$. Consider a finitely-generated
extension $F/E$ of transcendence degree $n$. By Lang's theorem
\cite{L1}, $E$ is of type $C _{1}$, whence, by the
Lang-Nagata-Tsen theorem \cite{Na}, $F$ is of type $C _{n+1}$. The
assumptions on $F$ and $E$ also imply the existence of a discrete
valuation $\omega ^{\prime }$ of $F$ extending $\omega $, such
that $\widehat F/\widehat E$ is a finitely-generated extension of
transcendence degree $n$ (when $F/E$ is purely transcendental, one
may take as $\omega ^{\prime }$ the restricted Gauss prolongation
of $\omega $ on $F$). Thus it follows that $[\widehat F ^{\prime
}\colon \widehat F ^{\prime p}] = p ^{n}$, for every finite
extension $F ^{\prime }/F$. This enables one to deduce (e.g., from
\cite{Ch4}, Lemmas~3.1 and 4.3) that if $(L, w)$ is a
Henselization of $(F, \omega ^{\prime })$, then abrd$_{p}(L) \le
{\rm Brd}_{p}(F)$. Hence, Conjecture \ref{conj7.7} and the $C
_{n+1}$ type of $F$ require that abrd$_{p}(L) \le n$. On the other
hand, $(L, w)/(F, \omega ^{\prime })$ is immediate, so $[\widehat
L\colon \widehat L ^{p}] = p ^{n}$, and by Theorem \ref{theo2.1}
(b), Brd$_{p}(L) \ge n$. Thus the assertion that Brd$_{p}(L) = n$
can be viewed as a special case of Conjecture \ref{conj7.7}.

\vskip0.38truecm \emph{Acknowledgment.} I would like to thank the
referee for the careful reading of an earlier version of this
paper, and for a number of suggestions used for improving the
organization (and other aspects) of its presentation. The paper
presents a research partially supported by Grant KP-06 N 32/1 of
07.12.2019 "Groups and Rings - Theory and Applications" of the
Bulgarian National Science Fund.

\end{document}